\documentclass{amsart}

\theoremstyle{definition}

\theoremstyle{remark}

\numberwithin{equation}{section}

\usepackage{graphicx,amsthm}
\usepackage{verbatim,amsmath,amscd,amssymb,color}
\usepackage[mathscr]{eucal}
\input xy
\xyoption{all}
\begin{document}
\renewcommand{\labelenumi}{$($\roman{enumi}$)$}
\renewcommand{\labelenumii}{$(${\rm \alph{enumii}}$)$}
\font\germ=eufm10
\newcommand{\cI}{{\mathcal I}}
\newcommand{\cA}{{\mathcal A}}
\newcommand{\cB}{{\mathcal B}}
\newcommand{\cC}{{\mathcal C}}
\newcommand{\cD}{{\mathcal D}}
\newcommand{\cE}{{\mathcal E}}
\newcommand{\cF}{{\mathcal F}}
\newcommand{\cG}{{\mathcal G}}
\newcommand{\cH}{{\mathcal H}}
\newcommand{\cK}{{\mathcal K}}
\newcommand{\cL}{{\mathcal L}}
\newcommand{\cM}{{\mathcal M}}
\newcommand{\cN}{{\mathcal N}}
\newcommand{\cO}{{\mathcal O}}
\newcommand{\cR}{{\mathcal R}}
\newcommand{\cS}{{\mathcal S}}
\newcommand{\cV}{{\mathcal V}}
\newcommand{\cX}{{\mathcal X}}
\newcommand{\fra}{\mathfrak a}
\newcommand{\frb}{\mathfrak b}
\newcommand{\frc}{\mathfrak c}
\newcommand{\frd}{\mathfrak d}
\newcommand{\fre}{\mathfrak e}
\newcommand{\frf}{\mathfrak f}
\newcommand{\frg}{\mathfrak g}
\newcommand{\frh}{\mathfrak h}
\newcommand{\fri}{\mathfrak i}
\newcommand{\frj}{\mathfrak j}
\newcommand{\frk}{\mathfrak k}
\newcommand{\frI}{\mathfrak I}
\newcommand{\fm}{\mathfrak m}
\newcommand{\frn}{\mathfrak n}
\newcommand{\frp}{\mathfrak p}
\newcommand{\fq}{\mathfrak q}
\newcommand{\frr}{\mathfrak r}
\newcommand{\frs}{\mathfrak s}
\newcommand{\frt}{\mathfrak t}
\newcommand{\fru}{\mathfrak u}
\newcommand{\frA}{\mathfrak A}
\newcommand{\frB}{\mathfrak B}
\newcommand{\frF}{\mathfrak F}
\newcommand{\frG}{\mathfrak G}
\newcommand{\frH}{\mathfrak H}
\newcommand{\frJ}{\mathfrak J}
\newcommand{\frN}{\mathfrak N}
\newcommand{\frP}{\mathfrak P}
\newcommand{\frT}{\mathfrak T}
\newcommand{\frU}{\mathfrak U}
\newcommand{\frV}{\mathfrak V}
\newcommand{\frX}{\mathfrak X}
\newcommand{\frY}{\mathfrak Y}
\newcommand{\frZ}{\mathfrak Z}
\newcommand{\rA}{\mathrm{A}}
\newcommand{\rC}{\mathrm{C}}
\newcommand{\rd}{\mathrm{d}}
\newcommand{\rB}{\mathrm{B}}
\newcommand{\rD}{\mathrm{D}}
\newcommand{\rE}{\mathrm{E}}
\newcommand{\rH}{\mathrm{H}}
\newcommand{\rK}{\mathrm{K}}
\newcommand{\rL}{\mathrm{L}}
\newcommand{\rM}{\mathrm{M}}
\newcommand{\rN}{\mathrm{N}}
\newcommand{\rR}{\mathrm{R}}
\newcommand{\rT}{\mathrm{T}}
\newcommand{\rZ}{\mathrm{Z}}
\newcommand{\bbA}{\mathbb A}
\newcommand{\bbB}{\mathbb B}
\newcommand{\bbC}{\mathbb C}
\newcommand{\bbG}{\mathbb G}
\newcommand{\bbF}{\mathbb F}
\newcommand{\bbH}{\mathbb H}
\newcommand{\bbP}{\mathbb P}
\newcommand{\bbN}{\mathbb N}
\newcommand{\bbQ}{\mathbb Q}
\newcommand{\bbR}{\mathbb R}
\newcommand{\bbV}{\mathbb V}
\newcommand{\bbZ}{\mathbb Z}
\newcommand{\adj}{\operatorname{adj}}
\newcommand{\Ad}{\mathrm{Ad}}
\newcommand{\Ann}{\mathrm{Ann}}
\newcommand{\rcris}{\mathrm{cris}}
\newcommand{\ch}{\mathrm{ch}}
\newcommand{\coker}{\mathrm{coker}}
\newcommand{\diag}{\mathrm{diag}}
\newcommand{\Diff}{\mathrm{Diff}}
\newcommand{\Dist}{\mathrm{Dist}}
\newcommand{\rDR}{\mathrm{DR}}
\newcommand{\ev}{\mathrm{ev}}
\newcommand{\Ext}{\mathrm{Ext}}
\newcommand{\cExt}{\mathcal{E}xt}
\newcommand{\fin}{\mathrm{fin}}
\newcommand{\Frac}{\mathrm{Frac}}
\newcommand{\GL}{\mathrm{GL}}
\newcommand{\Hom}{\mathrm{Hom}}
\newcommand{\hd}{\mathrm{hd}}
\newcommand{\rht}{\mathrm{ht}}
\newcommand{\id}{\mathrm{id}}
\newcommand{\im}{\mathrm{im}}
\newcommand{\inc}{\mathrm{inc}}
\newcommand{\ind}{\mathrm{ind}}
\newcommand{\coind}{\mathrm{coind}}
\newcommand{\Lie}{\mathrm{Lie}}
\newcommand{\Max}{\mathrm{Max}}
\newcommand{\mult}{\mathrm{mult}}
\newcommand{\op}{\mathrm{op}}
\newcommand{\ord}{\mathrm{ord}}
\newcommand{\pt}{\mathrm{pt}}
\newcommand{\qt}{\mathrm{qt}}
\newcommand{\rad}{\mathrm{rad}}
\newcommand{\res}{\mathrm{res}}
\newcommand{\rgt}{\mathrm{rgt}}
\newcommand{\rk}{\mathrm{rk}}
\newcommand{\SL}{\mathrm{SL}}
\newcommand{\soc}{\mathrm{soc}}
\newcommand{\Spec}{\mathrm{Spec}}
\newcommand{\St}{\mathrm{St}}
\newcommand{\supp}{\mathrm{supp}}
\newcommand{\Tor}{\mathrm{Tor}}
\newcommand{\Tr}{\mathrm{Tr}}
\newcommand{\wt}{\mathrm{wt}}
\newcommand{\Ab}{\mathbf{Ab}}
\newcommand{\Alg}{\mathbf{Alg}}
\newcommand{\Grp}{\mathbf{Grp}}
\newcommand{\Mod}{\mathbf{Mod}}
\newcommand{\Sch}{\mathbf{Sch}}\newcommand{\bfmod}{{\bf mod}}
\newcommand{\Qc}{\mathbf{Qc}}
\newcommand{\Rng}{\mathbf{Rng}}
\newcommand{\Top}{\mathbf{Top}}
\newcommand{\Var}{\mathbf{Var}}
\newcommand{\gromega}{\langle\omega\rangle}
\newcommand{\lbr}{\begin{bmatrix}}
\newcommand{\rbr}{\end{bmatrix}}
\newcommand{\cd}{commutative diagram }
\newcommand{\SpS}{spectral sequence}
\newcommand\C{\mathbb C}
\newcommand\hh{{\hat{H}}}
\newcommand\eh{{\hat{E}}}
\newcommand\F{\mathbb F}
\newcommand\fh{{\hat{F}}}
\newcommand\Z{{\mathbb Z}}
\newcommand\Zn{\Z_{\geq0}}
\newcommand\et[1]{\tilde{e}_{#1}}
\newcommand\ft[1]{\tilde{f}_{#1}}

\def\ge{\frg}
\def\AA{{\mathcal A}}
\def\al{\alpha}
\def\bq{B_q(\ge)}
\def\bqm{B_q^-(\ge)}
\def\bqz{B_q^0(\ge)}
\def\bqp{B_q^+(\ge)}
\def\beneme{\begin{enumerate}}
\def\beq{\begin{equation}}
\def\beqn{\begin{eqnarray}}
\def\beqnn{\begin{eqnarray*}}
\def\bigsl{{\hbox{\fontD \char'54}}}
\def\bbra#1,#2,#3{\left\{\begin{array}{c}\hspace{-5pt}
#1;#2\\ \hspace{-5pt}#3\end{array}\hspace{-5pt}\right\}}
\def\cd{\cdots}
\def\CC{\mathbb{C}}
\def\CBL{\cB_L(\TY(B,1,n+1))}
\def\CBM{\cB_M(\TY(B,1,n+1))}
\def\CVL{\cV_L(\TY(D,1,n+1))}
\def\CVM{\cV_M(\TY(D,1,n+1))}
\def\ddd{\hbox{\germ D}}
\def\del{\delta}
\def\Del{\Delta}
\def\Delr{\Delta^{(r)}}
\def\Dell{\Delta^{(l)}}
\def\Delb{\Delta^{(b)}}
\def\Deli{\Delta^{(i)}}
\def\Delre{\Delta^{\rm re}}
\def\ei{e_i}
\def\eit{\tilde{e}_i}
\def\eneme{\end{enumerate}}
\def\ep{\epsilon}
\def\eeq{\end{equation}}
\def\eeqn{\end{eqnarray}}
\def\eeqnn{\end{eqnarray*}}
\def\fit{\tilde{f}_i}
\def\FF{{\rm F}}
\def\ft{\tilde{f}_}
\def\gau#1,#2{\left[\begin{array}{c}\hspace{-5pt}#1\\
\hspace{-5pt}#2\end{array}\hspace{-5pt}\right]}
\def\gl{\hbox{\germ gl}}
\def\hom{{\hbox{Hom}}}
\def\ify{\infty}
\def\io{\iota}
\def\kp{k^{(+)}}
\def\km{k^{(-)}}
\def\llra{\relbar\joinrel\relbar\joinrel\relbar\joinrel\rightarrow}
\def\lan{\langle}
\def\lar{\longrightarrow}
\def\max{{\rm max}}
\def\lm{\lambda}
\def\Lm{\Lambda}
\def\mapright#1{\smash{\mathop{\longrightarrow}\limits^{#1}}}
\def\Mapright#1{\smash{\mathop{\Longrightarrow}\limits^{#1}}}
\def\mm{{\bf{\rm m}}}
\def\nd{\noindent}
\def\nn{\nonumber}
\def\nnn{\hbox{\germ n}}
\def\catob{{\mathcal O}(B)}
\def\oint{{\mathcal O}_{\rm int}(\ge)}
\def\ot{\otimes}
\def\op{\oplus}
\def\opi{\ovl\pi_{\lm}}
\def\osigma{\ovl\sigma}
\def\ovl{\overline}
\def\plm{\Psi^{(\lm)}_{\io}}
\def\qq{\qquad}
\def\q{\quad}
\def\qed{\hfill\framebox[2mm]{}}
\def\QQ{\mathbb Q}
\def\qi{q_i}
\def\qii{q_i^{-1}}
\def\ra{\rightarrow}
\def\ran{\rangle}
\def\rlm{r_{\lm}}
\def\ssl{\hbox{\germ sl}}
\def\slh{\widehat{\ssl_2}}
\def\ti{t_i}
\def\tii{t_i^{-1}}
\def\til{\tilde}
\def\tm{\times}
\def\tt{\frt}
\def\TY(#1,#2,#3){#1^{(#2)}_{#3}}
\def\ua{U_{\AA}}
\def\ue{U_{\vep}}
\def\uq{U_q(\ge)}
\def\uqp{U'_q(\ge)}
\def\ufin{U^{\rm fin}_{\vep}}
\def\ufinp{(U^{\rm fin}_{\vep})^+}
\def\ufinm{(U^{\rm fin}_{\vep})^-}
\def\ufinz{(U^{\rm fin}_{\vep})^0}
\def\uqm{U^-_q(\ge)}
\def\uqmq{{U^-_q(\ge)}_{\bf Q}}
\def\uqpm{U^{\pm}_q(\ge)}
\def\uqq{U_{\bf Q}^-(\ge)}
\def\uqz{U^-_{\bf Z}(\ge)}
\def\ures{U^{\rm res}_{\AA}}
\def\urese{U^{\rm res}_{\vep}}
\def\uresez{U^{\rm res}_{\vep,\ZZ}}
\def\util{\widetilde\uq}
\def\uup{U^{\geq}}
\def\ulow{U^{\leq}}
\def\bup{B^{\geq}}
\def\blow{\ovl B^{\leq}}
\def\vep{\varepsilon}
\def\vp{\varphi}
\def\vpi{\varphi^{-1}}
\def\VV{{\mathcal V}}
\def\xii{\xi^{(i)}}
\def\Xiioi{\Xi_{\io}^{(i)}}
\def\W1{W(\varpi_1)}
\def\WW{{\mathcal W}}
\def\wt{{\rm wt}}
\def\wtil{\widetilde}
\def\what{\widehat}
\def\wpi{\widehat\pi_{\lm}}
\def\ZZ{\mathbb Z}
\def\RR{\mathbb R}

\def\m@th{\mathsurround=0pt}
\def\fsquare(#1,#2){
\hbox{\vrule$\hskip-0.4pt\vcenter to #1{\normalbaselines\m@th
\hrule\vfil\hbox to #1{\hfill$\scriptstyle #2$\hfill}\vfil\hrule}$\hskip-0.4pt
\vrule}}

\newtheorem{thm}{Theorem}[section]
\newtheorem{pro}[thm]{Proposition}
\newtheorem{lem}[thm]{Lemma}
\newtheorem{ex}[thm]{Example}
\newtheorem{cor}[thm]{Corollary}
\newtheorem{conj}[thm]{Conjecture}
\theoremstyle{definition}
\newtheorem{df}[thm]{Definition}

\newcommand{\cmt}{\marginpar}
\newcommand{\seteq}{\mathbin{:=}}
\newcommand{\cl}{\colon}
\newcommand{\be}{\begin{enumerate}}
\newcommand{\ee}{\end{enumerate}}
\newcommand{\bnum}{\be[{\rm (i)}]}
\newcommand{\enum}{\ee}
\newcommand{\ro}{{\rm(}}
\newcommand{\rf}{{\rm)}}
\newcommand{\set}[2]{\left\{#1\,\vert\,#2\right\}}
\newcommand{\sbigoplus}{{\mbox{\small{$\bigoplus$}}}}
\newcommand{\ba}{\begin{array}}
\newcommand{\ea}{\end{array}}
\newcommand{\on}{\operatorname}
\newcommand{\eq}{\begin{eqnarray}}
\newcommand{\eneq}{\end{eqnarray}}
\newcommand{\hs}{\hspace*}

\title[Ultra-discretization 
of the $\TY(G,1,2)$-Geometric Crystals]
{Ultra-discretization 
of the $\TY(G,1,2)$-Geometric Crystals to the
$\TY(D,3,4)$-Perfect Crystals}


\author{Toshiki N\textsc{akashima}}
\address{Department of Mathematics, 
Sophia University, Kioicho 7-1, Chiyoda-ku, Tokyo 102-8554,
Japan}
\email{toshiki@mm.sophia.ac.jp}
\thanks{supported in part by JSPS Grants 
in Aid for Scientific Research $\sharp 19540050$.}

\subjclass{Primary 17B37; 17B67; Secondary 22E65; 14M15}
\date{}


\keywords{geometric crystal, perfect crystal, 
ultra-discretization. }

\begin{abstract}
We obtain the affirmative answer 
to the conjecture 
in \cite{N3}. More, precisely, 
let $\chi\seteq(\cV,\{e_i\},\{\gamma_i\},
\{\vep_i\})$
be the affine geometric crystal of type 
$\TY(G,1,2)$
in \cite{N2} and 
${\mathcal UD}(\chi,T,\theta)$ 
a ultra-discretization 
of $\chi$ with respect to a 
certain positive structure $\theta$.
Then we show that 
${\mathcal UD}(\chi,T,\theta)$ 
is isomorphic to the limit of 
coherent family of perfect crystals of 
type $\TY(D,3,4)$ in \cite{KMOY}.
\end{abstract}

\maketitle
\renewcommand{\thesection}{\arabic{section}}
\section{Introduction}
\setcounter{equation}{0}
\renewcommand{\theequation}{\thesection.\arabic{equation}}

In \cite{KMN1}, 
we introduced the notion of perfect crystal, 
which holds several nice properties, 
{\it e.g.,}, the existence of the isomorphism of 
crystals:
\[
 B(\lm)\cong B(\sigma(\lm))\ot B,
\]
where $B$ is a perfect crystal of level 
$l\in\bbZ_{>0}$,
$B(\lm)$ is the crystal of the integrable 
highest weight module of a quantum affine group with
the level $l$ highest weight $\lm$ 
and $\sigma$ is a certain bijection on dominant
weights. Iterating this isomorphism, one can 
get the so-called Kyoto path model for $B(\lm)$,
which plays an crucial role in 
calculating the one-point functions for
vertex-type lattice models 
(\cite{KMN1},\cite{KMN2}).

In \cite{KMN2} perfect crystals with 
arbitrary level has been constructed explicitly
for affine Kac-Moody algebra of type
$\TY(A,1,n)$, $\TY(B,1,n)$, $\TY(C,1,n)$, 
$\TY(D,1,n)$, $\TY(D,2,n+1)$, $\TY(A,2,2n-1)$
and $\TY(A,2,2n)$. In \cite{Y}, the $\TY(G,1,2)$
case has been accomplished. But, so far 
the other cases except $\TY(D,3,4)$ have 
not yet been obtained. In the recent work \cite{KMOY}, 
they constructed the perfect crystal of type 
$\TY(D,3,4)$ with arbitrary level explicitly.
A coherent family of perfect crystals is defined in 
\cite{KKM} and it has been shown that 
the perfect crystals 
in \cite{KMN2} constitute a coherent family.
A coherent family $\{B_l\}_{l\geq1}$ 
of perfect crystals $B_l$ possesses a limit
$B_\ify$ which still keeps a structure of crystal.
This has a similar property to $B_l$, that is, 
there exists the isomorphism of crystals:
\[
 B(\ify)\cong B(\ify)\ot B_\ify,
\]
where $B(\ify)$ is the crystal of the 
nilpotent subalgebra
$\uqm$ of a quantum affine algebra $\uq$. 
An iteration of the isomorphism also produces 
a path model of $B(\ify)$(\cite{KKM}).
It is shown in \cite{KMOY} 
that the obtained perfect crystals consists of a 
coherent family and the structure of the limit 
$B_\ify$ has been described explicitly.

Geometric crystal is 
 an object defined over certain algebraic
(or ind-)variety which seems to be a kind of 
geometric lifting of Kashiwara's crystal. 
It is defined
in \cite{BK} for reductive algebraic groups and 
is extended to general Kac-Moody cases in \cite{N}. 
For a fixed Cartan data
$(A,\{\al_i\}_{i\in I},\{\al^\vee_i\}_{\i\in I})$, 
a geometric crystal consists of
an ind-variety $X$ over 
the complex number $\bbC$, 
a rational $\bbC^\times$-action 
$e_i:\bbC^\times\times X\longrightarrow X$ and 
rational functions 
$\gamma_i,\vep_i:X\longrightarrow 
\bbC$ $(i\in I)$,
which satisfy the conditions as in 
Definition \ref{def-gc}.
It has many similarity to 
the theory of crystals, {\it e.g.,}
some product structure, 
Weyl group actions, R-matrices, 
{\it etc}. Moreover,
one has a direct connection
between geometric crystals and free crystals, 
called tropicalization/
ultra-discretization procedure (see \S2).
Here let us explain this procedure.
For an algebraic torus $T'$ and a birational 
morphism $\theta:T'\to X$, the pair $(T',\theta)$
is positive if it satisfies 
the conditions as in Sect.2, roughly speaking:
Through the morphism $\theta$, 
we can induce a geometric crystal structure on $T'$ 
from $X$ and 
express the data $e_i^c$, $\gamma_i$ and $\vep_i$
$(i\in I)$
using the coordinate of $T'$ explicitly.
In case each of them are expressed as a ratio of 
positive polynomials, it is said that 
$(T',\theta)$ is a positive structure of the 
geometric crystal $(X,\{e_i\},\{\gamma_i\},
\{\vep_i\})$.
Then by using a map $v:\bbC(c)\setminus\{0\}\to \bbZ$ 
($v(f)\seteq {\rm deg}(f)$), we can define a 
morphism $T'\to \bbZ^m$ ($m=\dim T'=\dim X$), which 
defines the so-called ultra-discretization functor. 
If $\theta:T'\to X$ is a positive structure on $X$, 
then we obtain a Kashiwara's crystal from $X$ by 
applying the ultra-discretization functor(\cite{BK}).

Let $G$ (resp.~$\ge=\lan \tt,e_i,f_i\ran_{i\in I}$) be the affine 
Kac-Moody group 
(resp.~algebra) associated with a 
generalized Cartan matrix $A=(a_{ij})_{i,j\in I}$. 
Let $B^\pm$ be fixed Borel subgroups and $T$ the 
maximal torus such that $B^+\cap B^-=T$. 
Set $y_i(c)\seteq\exp(cf_i)$,
and let $\al_i^\vee(c)\in T$ 
be the image of $c\in\bbC^\times$
by the group morphism $\bbC^\times\to T$ induced by
the simple coroot $\alpha_i^\vee$ as in \ref{KM}.
We set $Y_i(c)\seteq y_i(c^{-1})\,\al_i^\vee(c)=\al_i^\vee(c)\,y_i(c)$.
Let $W$ (resp.~$\wtil W$) be the Weyl group 
(resp.~the extended Weyl group) associated with $\ge$.
The Schubert cell $X_w\seteq BwB/B$ $(w=s_{i_1}\cd s_{i_k}\in W)$ 
is birationally isomorphic to the variety
\[
 B^-_w\seteq\set{Y_{i_1}(x_1)\cd Y_{i_k}(x_k)}%
{x_1,\cd,x_k\in \bbC^\times}\subset B^-,
\]
and $X_w$ has a natural geometric crystal structure
(\cite{BK}, \cite{N}).

We choose $0\in I$ as in \cite{K0},
and let $\{\varpi_i\}_{i\in I\setminus\{0\}}$
be the set of level $0$ fundamental weights.
Let $W(\varpi_i)$ be the fundamental representation
of $U_q(\ge)$ with $\varpi_i$ as an extremal weight (\cite{K0}).
Let us denote its specialization at $q=1$
by the same notation $W(\varpi_i)$. 
It is a finite-dimensional $\ge$-module.
Let $\bbP(\varpi_i)$ be the projective space
$(W(\varpi_i)\setminus\{0\})/\bbC^\times$.

For any $i\in I$, define 
$c_i^\vee\seteq
\mathrm{max }(1,\frac{2}{(\al_i,\al_i)}).$
Then the translation $t(c^\vee_i\varpi_i)$ 
belongs to $\widetilde W$ (see \cite{KNO}).
For a subset $J$ of $I$, let us denote by
$\ge_J$ the subalgebra of $\ge$ generated by
$\{e_i,f_i\}_{i\in J}$.
For an integral weight $\mu$, define 
$I(\mu)\seteq\set{j\in I}{\lan \al^\vee_j,\mu\ran\geq0}$. 

Here we state the conjecture given in \cite{KNO}:
\begin{conj}[\cite{KNO}]
For any $i\in I$, 
there exist a unique variety $X$ endowed with
a positive $\ge$-geometric crystal structure and a
rational mapping $\pi\cl X\longrightarrow
\bbP(\varpi_i)$ satisfying the following property:
\begin{enumerate}
\item
for an arbitrary extremal vector $u\in W(\varpi_i)_\mu$,
writing the translation
$t(c_i^\vee\mu)$ as $\io w\in 
\wtil W$ with a Dynkin diagram automorphism $\io$
and $w=s_{i_1}\cd s_{i_k}$,
there exists a birational mapping
$\xi\cl B^-_w\longrightarrow X$
such that $\xi$ is a morphism of $\ge_{I(\mu)}$-geometric crystals
and that
the composition
$\pi\circ\xi\cl B^-_w\to \bbP(\varpi_i)$
coincides with
$Y_{i_1}(x_1)\cdots Y_{i_k}(x_k)\mapsto 
Y_{i_1}(x_1)\cd Y_{i_k}(x_k)\ovl u$,
where $\ovl u$ is the line including $u$,
\item
the ultra-discretization(see Sect.2) of $X$
is isomorphic to the crystal $B_\infty(\varpi_i)$
of the Langlands dual $\ge^L$.
\end{enumerate}
\end{conj}
In \cite{KNO}, the cases $i=1$ and 
$\ge=\TY(A,1,n), \TY(B,1,n),\TY(C,1,n),\TY(D,1,n),$
$\TY(A,2,2n-1),\TY(A,2,2n),\TY(D,2,n+1)$ have been 
resolved, that is, 
certain positive geometric crystal 
$\cV(\ge)$ associated 
with the fundamental representation $W(\varpi_1)$ 
for the above affine Lie algebras has been 
constructed and it was shown 
that the ultra-discretization 
limit of $\cV(\ge)$ 
is isomorphic to the limit of 
the coherent family 
of perfect crystals as above
for $\ge^L$ the Langlands dual of $\ge$.
In \cite{N3} for
the case $i=1$ and $\ge=\TY(G,1,2)$, a 
positive geometric
crystal $\cV$ was constructed. However, 
the ultra-discretization of the 
geometric crystal has not been given there,
though it was conjectured that 
the ultra-discretization of $\cV$ is isomorphic 
to $B_\ify$ as in \cite{KMOY}.

In this article, we shall describe the structure of 
the crystal obtained
by ultra-discretization process from the 
geometric crystal $\cV$ 
for $\ge=\TY(G,1,2)$ in \cite{N3}.
Finally, we shall show that the crystal is isomorphic
to $B_\ify$ as in \cite{KMOY}.
\renewcommand{\thesection}{\arabic{section}}
\section{Geometric crystals}
\setcounter{equation}{0}
\renewcommand{\theequation}{\thesection.\arabic{equation}}

In this section, 
we review Kac-Moody groups and geometric crystals
following 
\cite{PK}, \cite{Ku2}, \cite{BK}
\subsection{Kac-Moody algebras and Kac-Moody groups}
\label{KM}
Fix a symmetrizable generalized Cartan matrix
 $A=(a_{ij})_{i,j\in I}$ with a finite index set $I$.
Let $(\tt,\{\al_i\}_{i\in I},\{\al^\vee_i\}_{i\in I})$ 
be the associated
root data, where ${\tt}$ is a vector space 
over $\bbC$ and
$\{\al_i\}_{i\in I}\subset\tt^*$ and 
$\{\al^\vee_i\}_{i\in I}\subset\tt$
are linearly independent 
satisfying $\al_j(\al^\vee_i)=a_{ij}$.

The Kac-Moody Lie algebra $\ge=\ge(A)$ associated with $A$
is the Lie algebra over $\bbC$ generated by $\tt$, the 
Chevalley generators $e_i$ and $f_i$ $(i\in I)$
with the usual defining relations (\cite{KP},\cite{PK}).
There is the root space decomposition 
$\ge=\bigoplus_{\al\in \tt^*}\ge_{\al}$.
Denote the set of roots by 
$\Delta:=\{\al\in \tt^*|\al\ne0,\,\,\ge_{\al}\ne(0)\}$.
Set $Q=\sum_i\bbZ \al_i$, $Q_+=\sum_i\bbZ_{\geq0} \al_i$,
$Q^\vee:=\sum_i\bbZ \al^\vee_i$
and $\Delta_+:=\Delta\cap Q_+$.
An element of $\Delta_+$ is called 
a {\it positive root}.
Let $P\subset \tt^*$ be a weight lattice such that 
$\bbC\ot P=\tt^*$, whose element is called a
weight.

Define simple reflections $s_i\in{\rm Aut}(\tt)$ $(i\in I)$ by
$s_i(h):=h-\al_i(h)\al^\vee_i$, which generate the Weyl group $W$.
It induces the action of $W$ on $\tt^*$ by
$s_i(\lm):=\lm-\lm(\al^\vee_i)\al_i$.
Set $\Delre:=\{w(\al_i)|w\in W,\,\,i\in I\}$, whose element 
is called a real root.

Let $\ge'$ be the derived Lie algebra 
of $\ge$ and let 
$G$ be the Kac-Moody group associated 
with $\ge'$(\cite{PK}).
Let $U_{\al}:=\exp\ge_{\al}$ $(\al\in \Delre)$
be the one-parameter subgroup of $G$.
The group $G$ is generated by $U_{\al}$ $(\al\in \Delre)$.
Let $U^{\pm}$ be the subgroup generated by $U_{\pm\al}$
($\al\in \Delre_+=\Delre\cap Q_+$), {\it i.e.,}
$U^{\pm}:=\lan U_{\pm\al}|\al\in\Del^{\rm re}_+\ran$.

For any $i\in I$, there exists a unique homomorphism;
$\phi_i:SL_2(\bbC)\rightarrow G$ such that
\[
\hspace{-2pt}\phi_i\left(
\left(
\begin{array}{cc}
c&0\\
0&c^{-1}
\end{array}
\right)\right)=c^{\al^\vee_i},\,
\phi_i\left(
\left(
\begin{array}{cc}
1&t\\
0&1
\end{array}
\right)\right)=\exp(t e_i),\,
 \phi_i\left(
\left(
\begin{array}{cc}
1&0\\
t&1
\end{array}
\right)\right)=\exp(t f_i).
\]
where $c\in\bbC^\times$ and $t\in\bbC$.
Set $\al^\vee_i(c):=c^{\al^\vee_i}$,
$x_i(t):=\exp{(t e_i)}$, $y_i(t):=\exp{(t f_i)}$, 
$G_i:=\phi_i(SL_2(\bbC))$,
$T_i:=\phi_i(\{{\rm diag}(c,c^{-1})\vert 
c\in\bbC^{\vee}\})$ 
and 
$N_i:=N_{G_i}(T_i)$. Let
$T$ (resp. $N$) be the subgroup of $G$ 
with the Lie algebra $\tt$
(resp. generated by the $N_i$'s), 
which is called a {\it maximal torus} in $G$, and let
$B^{\pm}=U^{\pm}T$ be the Borel subgroup of $G$.
We have the isomorphism
$\phi:W\mapright{\sim}N/T$ defined by $\phi(s_i)=N_iT/T$.
An element $\ovl s_i:=x_i(-1)y_i(1)x_i(-1)
=\phi_i\left(
\left(
\begin{array}{cc}
0&\pm1\\
\mp1&0
\end{array}
\right)\right)$ is in 
$N_G(T)$, which is a representative of 
$s_i\in W=N_G(T)/T$. 

\subsection{Geometric crystals}
Let $W$ be the  Weyl group associated with $\ge$. 
Define $R(w)$ for $w\in W$ by
\[
 R(w):=\{(i_1,i_2,\cd,i_l)\in I^l|w=s_{i_1}s_{i_2}\cd s_{i_l}\},
\]
where $l$ is the length of $w$.
Then $R(w)$ is the set of reduced words of $w$.

Let $X$ be an ind-variety , 
{$\gamma_i:X\rightarrow \bbC$} and 
$\vep_i:X\longrightarrow \bbC$ ($i\in I$) 
rational functions on $X$, and
{$e_i:\bbC^\times \times X\longrightarrow X$}
$((c,x)\mapsto e^c_i(x))$ a
rational $\bbC^\times$-action.

For a word ${\bf i}=(i_1,\cd,i_l)\in R(w)$ 
$(w\in W)$, set 
$\al^{(j)}:=s_{i_l}\cd s_{i_{j+1}}(\al_{i_j})$ 
$(1\leq j\leq l)$ and 
\begin{eqnarray*}
e_{\bf i}:&T\times X\rightarrow &X\\
&(t,x)\mapsto &e_{\bf i}^t(x):=e_{i_1}^{\al^{(1)}(t)}
e_{i_2}^{\al^{(2)}(t)}\cd e_{i_l}^{\al^{(l)}(t)}(x).
\label{tx}
\end{eqnarray*}
\begin{df}
\label{def-gc}
A quadruple $(X,\{e_i\}_{i\in I},\{\gamma_i,\}_{i\in I},
\{\vep_i\}_{i\in I})$ is a 
$G$ (or $\ge$)-\\{\it geometric} {\it crystal} 
if
\begin{enumerate}
\item
$\{1\}\times X\subset dom(e_i)$ 
for any $i\in I$.
\item
$\gamma_j(e^c_i(x))=c^{a_{ij}}\gamma_j(x)$.
\item
{$e_{\bf i}=e_{\bf i'}$}
for any 
$w\in W$, ${\bf i}$.
${\bf i'}\in R(w)$.
\item
$\vep_i(e_i^c(x))=c^{-1}\vep_i(x)$.
\end{enumerate}
\end{df}
Note that the condition (iii) as above is 
equivalent to the following so-called 
{\it Verma relations}:
\[
 \begin{array}{lll}
&\hspace{-20pt}e^{c_1}_{i}e^{c_2}_{j}
=e^{c_2}_{j}e^{c_1}_{i}&
{\rm if }\,\,a_{ij}=a_{ji}=0,\\
&\hspace{-20pt} e^{c_1}_{i}e^{c_1c_2}_{j}e^{c_2}_{i}
=e^{c_2}_{j}e^{c_1c_2}_{i}e^{c_1}_{j}&
{\rm if }\,\,a_{ij}=a_{ji}=-1,\\
&\hspace{-20pt}
e^{c_1}_{i}e^{c^2_1c_2}_{j}e^{c_1c_2}_{i}e^{c_2}_{j}
=e^{c_2}_{j}e^{c_1c_2}_{i}e^{c^2_1c_2}_{j}e^{c_1}_{i}&
{\rm if }\,\,a_{ij}=-2,\,
a_{ji}=-1,\\
&\hspace{-20pt}
e^{c_1}_{i}e^{c^3_1c_2}_{j}e^{c^2_1c_2}_{i}
e^{c^3_1c^2_2}_{j}e^{c_1c_2}_{i}e^{c_2}_{j}
=e^{c_2}_{j}e^{c_1c_2}_{i}e^{c^3_1c^2_2}_{j}e^{c^2_1c_2}_{i}
e^{c^3_1c_2}_je^{c_1}_i&
{\rm if }\,\,a_{ij}=-3,\,
a_{ji}=-1,
\end{array}
\]
Note that the last formula is different from the one in 
\cite{BK}, \cite{N}, \cite{N2} which seems to be
incorrect. The formula here may be correct.

\subsection{Geometric crystal on Schubert cell}
\label{schubert}

Let $w\in W$ be a Weyl group element and take a 
reduced expression $w=s_{i_1}\cd s_{i_l}$. 
Let $X:=G/B$ be the flag
variety, which is an ind-variety 
and $X_w\subset X$ the
Schubert cell associated with $w$, which has 
a natural geometric crystal structure
(\cite{BK},\cite{N}).
For ${\bf i}:=(i_1,\cd,i_k)$, set 
\begin{equation}
B_{\bf i}^-
:=\{Y_{\bf i}(c_1,\cd,c_k)
:=Y_{i_1}(c_1)\cd Y_{i_l}(c_k)
\,\vert\, c_1\cd,c_k\in\bbC^\times\}\subset B^-,
\label{bw1}
\end{equation}
which has a geometric crystal structure(\cite{N})
isomorphic to $X_w$. 
The explicit forms of the action $e^c_i$, the rational 
function $\vep_i$  and $\gamma_i$ on 
$B_{\bf i}^-$ are given by
\begin{eqnarray}
&& e_i^c(Y_{i_1}(c_1)\cd Y_{i_l}(c_k))
=Y_{i_1}({\mathcal C}_1)\cd Y_{i_l}({\mathcal C}_k)),\nn \\
&&\text{where}\nn\\
&&{\mathcal C}_j:=
c_j\cdot \frac{\displaystyle \sum_{1\leq m\leq j,i_m=i}
 \frac{c}
{c_1^{a_{i_1,i}}\cd c_{m-1}^{a_{i_{m-1},i}}c_m}
+\sum_{j< m\leq k,i_m=i} \frac{1}
{c_1^{a_{i_1,i}}\cd c_{m-1}^{a_{i_{m-1},i}}c_m}}
{\displaystyle\sum_{1\leq m<j,i_m=i} 
 \frac{c}
{c_1^{a_{i_1,i}}\cd c_{m-1}^{a_{i_{m-1},i}}c_m}+
\mathop\sum_{j\leq m\leq k,i_m=i}  \frac{1}
{c_1^{a_{i_1,i}}\cd c_{m-1}^{a_{i_{m-1},i}}c_m}},
\label{eici}\\
&& \vep_i(Y_{i_1}(c_1)\cd Y_{i_l}(c_k))=
\sum_{1\leq m\leq k,i_m=i} \frac{1}
{c_1^{a_{i_1,i}}\cd c_{m-1}^{a_{i_{m-1},i}}c_m},
\label{vep-i}\\
&&\gamma_i(Y_{i_1}(c_1)\cd Y_{i_l}(c_k))
=c_1^{a_{i_1,i}}\cd c_k^{a_{i_k,i}}.
\label{gamma-i}
\end{eqnarray}

\subsection{Positive structure,\,\,
Ultra-discretizations \,\, and \,\,Tropicalizations}
\label{positive-str}

Let us recall the notions of 
positive structure, ultra-discretization and tropicalization.

The setting below is same as \cite{KNO}.
Let $T=(\bbC^\times)^l$ be an algebraic torus over $\bbC$ and 
$X^*(T):={\rm Hom}(T,\bbC^\times)\cong \ZZ^l$ 
(resp. $X_*(T):={\rm Hom}(\bbC^\times,T)\cong \ZZ^l$) 
be the lattice of characters
(resp. co-characters)
of $T$. 
Set $R:=\bbC(c)$ and define
$$
\begin{array}{cccc}
v:&R\setminus\{0\}&\longrightarrow &\ZZ\\
&f(c)&\mapsto
&{\rm deg}(f(c)),
\end{array}
$$
where $\rm deg$ is the degree of poles at $c=\ify$. 
Here note that for $f_1,f_2\in R\setminus\{0\}$, we have
\begin{equation}
v(f_1 f_2)=v(f_1)+v(f_2),\q
v\left(\frac{f_1}{f_2}\right)=v(f_1)-v(f_2)
\label{ff=f+f}
\end{equation}
A non-zero rational function on
an algebraic torus $T$ is called {\em positive} if
it is written as $g/h$ where
$g$ and $h$ are a positive linear combination of
characters of $T$.
\begin{df}
Let 
$f\cl T\rightarrow T'$ be 
a rational morphism between
two algebraic tori $T$ and 
$T'$.
We say that $f$ is {\em positive},
if $\chi\circ f$ is positive
for any character $\chi\cl T'\to \C$.
\end{df}
Denote by ${\rm Mor}^+(T,T')$ the set of 
positive rational morphisms from $T$ to $T'$.

\begin{lem}[\cite{BK}]
\label{TTT}
For any $f\in {\rm Mor}^+(T_1,T_2)$             
and $g\in {\rm Mor}^+(T_2,T_3)$, 
the composition $g\circ f$
is well-defined and belongs to ${\rm Mor}^+(T_1,T_3)$.
\end{lem}

By Lemma \ref{TTT}, we can define a category ${\mathcal T}_+$
whose objects are algebraic tori over $\bbC$ and arrows
are positive rational morphisms.

Let $f\cl T\rightarrow T'$ be a 
positive rational morphism
of algebraic tori $T$ and 
$T'$.
We define a map $\what f\cl X_*(T)\rightarrow X_*(T')$ by 
\[
\langle\chi,\what f(\xi)\rangle
=v(\chi\circ f\circ \xi),
\]
where $\chi\in X^*(T')$ and $\xi\in X_*(T)$.
\begin{lem}[\cite{BK}]
For any algebraic tori $T_1$, $T_2$, $T_3$, 
and positive rational morphisms 
$f\in {\rm Mor}^+(T_1,T_2)$, 
$g\in {\rm Mor}^+(T_2,T_3)$, we have
$\what{g\circ f}=\what g\circ\what f.$
\end{lem}
By this lemma, we obtain a functor 
\[
\begin{array}{cccc}
{\mathcal UD}:&{\mathcal T}_+&\longrightarrow &{{\hbox{\germ Set}}}\\
&T&\mapsto& X_*(T)\\
&(f:T\rightarrow T')&\mapsto& 
(\what f:X_*(T)\rightarrow X_*(T')))
\end{array}
\]

\begin{df}[\cite{BK}]
Let $\chi=(X,\{e_i\}_{i\in I},\{{\rm wt}_i\}_{i\in I},
\{\vep_i\}_{i\in I})$ be a 
geometric crystal, $T'$ an algebraic torus
and $\theta:T'\rightarrow X$ 
a birational isomorphism.
The isomorphism $\theta$ is called 
{\it positive structure} on
$\chi$ if it satisfies
\begin{enumerate}
\item for any $i\in I$ the rational functions
$\gamma_i\circ \theta:T'\rightarrow \bbC$ and 
$\vep_i\circ \theta:T'\rightarrow \bbC$ 
are positive.
\item
For any $i\in I$, the rational morphism 
$e_{i,\theta}:\bbC^\tm \tm T'\rightarrow T'$ defined by
$e_{i,\theta}(c,t)
:=\theta^{-1}\circ e_i^c\circ \theta(t)$
is positive.
\end{enumerate}
\end{df}
Let $\theta:T\rightarrow X$ be a positive structure on 
a geometric crystal $\chi=(X,\{e_i\}_{i\in I},$
$\{{\rm wt}_i\}_{i\in I},
\{\vep_i\}_{i\in I})$.
Applying the functor ${\mathcal UD}$ 
to positive rational morphisms
$e_{i,\theta}:\bbC^\tm \tm T'\rightarrow T'$ and
$\gamma\circ \theta:T'\ra T$
(the notations are
as above), we obtain
\begin{eqnarray*}
\til e_i&:=&{\mathcal UD}(e_{i,\theta}):
\ZZ\tm X_*(T) \rightarrow X_*(T)\\
{\rm wt}_i&:=&{\mathcal UD}(\gamma_i\circ\theta):
X_*(T')\rightarrow \bbZ,\\
\vep_i&:=&{\mathcal UD}(\vep_i\circ\theta):
X_*(T')\rightarrow \bbZ.
\end{eqnarray*}
Now, for given positive structure $\theta:T'\rightarrow X$
on a geometric crystal 
$\chi=(X,\{e_i\}_{i\in I},$
$\{{\rm wt}_i\}_{i\in I},
\{\vep_i\}_{i\in I})$, we associate 
the quadruple $(X_*(T'),\{\til e_i\}_{i\in I},
\{{\rm wt}_i\}_{i\in I},\{\vep_i\}_{i\in I})$
with a free pre-crystal structure (see \cite[2.2]{BK}) 
and denote it by ${\mathcal UD}_{\theta,T'}(\chi)$.
We have the following theorem:

\begin{thm}[\cite{BK}\cite{N}]
For any geometric crystal 
$\chi=(X,\{e_i\}_{i\in I},\{\gamma_i\}_{i\in I},$
$\{\vep_i\}_{i\in I})$ and positive structure
$\theta:T'\rightarrow X$, the associated pre-crystal
${\mathcal UD}_{\theta,T'}(\chi)=$\\
$(X_*(T'),\{e_i\}_{i\in I},\{{\rm wt}_i\}_{i\in I},
\{\vep_i\}_{i\in I})$ 
is a crystal {\rm (see \cite[2.2]{BK})}
\end{thm}

Now, let ${\mathcal GC}^+$ be a category whose 
object is a triplet
$(\chi,T',\theta)$ where 
$\chi=(X,\{e_i\},\{\gamma_i\},\{\vep_i\})$ 
is a geometric crystal and $\theta:T'\rightarrow X$ 
is a positive structure on $\chi$, and morphism
$f:(\chi_1,T'_1,\theta_1)\longrightarrow 
(\chi_2,T'_2,\theta_2)$ is given by a morphism 
$\vp:X_1\longrightarrow X_2$  
($\chi_i=(X_i,\cd)$) such that 
\[
f:=\theta_2^{-1}\circ\vp\circ\theta_1:T'_1\longrightarrow T'_2,
\]
is a positive rational morphism. Let ${\mathcal CR}$
be a category of crystals. 
Then by the theorem above, we have
\begin{cor}
\label{cor-posi}
$\mathcal UD_{\theta,T'}$ as above defines a functor
\begin{eqnarray*}
 {\mathcal UD}&:&{\mathcal GC}^+\longrightarrow {\mathcal CR},\\
&&(\chi,T',\theta)\mapsto X_*(T'),\\
&&(f:(\chi_1,T'_1,\theta_1)\rightarrow 
(\chi_2,T'_2,\theta_2))\mapsto
(\what f:X_*(T'_1)\rightarrow X_*(T'_2)).
\end{eqnarray*}

\end{cor}
We call the functor $\mathcal UD$
{\it ``ultra-discretization''} as \cite{N},\cite{N2}
instead of ``tropicalization'' as in \cite{BK}.
And 
for a crystal $B$, if there
exists a geometric crystal $\chi$ and a positive 
structure $\theta:T'\rightarrow X$ on $\chi$ such that 
${\mathcal UD}(\chi,T',\theta)\cong B$ as crystals, 
we call an object $(\chi,T',\theta)$ in ${\mathcal GC}^+$
a {\it tropicalization} of $B$, where 
it is not known that this correspondence is a functor.

\renewcommand{\thesection}{\arabic{section}}
\section{Limit of perfect crystals}
\label{limit}
\setcounter{equation}{0}
\renewcommand{\theequation}{\thesection.\arabic{equation}}
We review limit of perfect crystals following \cite{KKM}.
(See also \cite{KMN1},\cite{KMN2}).

\subsection{Crystals}

First we review the theory of crystals,
which is the notion obtained by
abstracting the combinatorial 
properties of crystal bases.
Let 
$(A,\{\al_i\}_{i\in I},\{\al^\vee_i\}_{\i\in I})$ be a 
Cartan data. 
\begin{df}
A {\it crystal} $B$ is a set endowed with the following maps:
\begin{eqnarray*}
&& {\rm wt}:B\lar P,\\
&&\vep_i:B\lar\ZZ\sqcup\{-\infty\},\q
  \vp_i:B\lar\ZZ\sqcup\{-\infty\} \q{\hbox{for}}\q i\in I,\\
&&\eit:B\sqcup\{0\}\lar B\sqcup\{0\},
\q\fit:B\sqcup\{0\}\lar B\sqcup\{0\}\q{\hbox{for}}\q i\in I,\\
&&\eit(0)=\fit(0)=0.
\end{eqnarray*}
those maps satisfy the following axioms: for
 all $b,b_1,b_2 \in B$, we have
\begin{eqnarray*}
&&\vp_i(b)=\vep_i(b)+\lan \al^\vee_i,{\rm wt}
(b)\ran,\\
&&\wt(\eit b)=\wt(b)+\al_i{\hbox{ if  }}\eit b\in B,\\
&&\wt(\fit b)=\wt(b)-\al_i{\hbox{ if  }}\fit b\in B,\\
&&\eit b_2=b_1 \Longleftrightarrow \fit b_1=b_2\,\,(\,b_1,b_2 \in B),\\
&&\vep_i(b)=-\ify
   \Longrightarrow \eit b=\fit b=0.
\end{eqnarray*}
\end{df}
The following tensor product structure 
is one of the most crucial properties of crystals.
\begin{thm}
\label{tensor}
Let $B_1$ and $B_2$ be crystals.
Set
$B_1\ot B_2:=
\{b_1\otimes b_2;\;b_j\in B_j\;(j=1,2)\}$. Then we have 
\begin{enumerate}
\item $B_1\ot B_2$ is a crystal.
\item
For $b_1\in B_1$ and $b_2\in B_2$, we have
$$
\tilde f_i(b_1\otimes b_2)=
\left\{\begin{array}{ll}\tilde f_ib_1\otimes b_2&
{\rm if}\;\varphi_i(b_1)>\vep_i(b_2),\\
b_1\otimes\tilde f_ib_2&{\rm if}\;
\varphi_i(b_1)\leq\vep_i(b_2).
\end{array}\right.
$$
$$
\tilde e_i(b_1\otimes b_2)=\left\{\begin{array}{ll}
b_1\otimes \tilde e_ib_2&
{\rm if}\;\varphi_i(b_1)<\vep_i(b_2),\\
\tilde e_ib_1\otimes b_2
&{\rm if}\;\varphi_i(b_1)\geq\vep_i(b_2),
\end{array}\right.
$$
\end{enumerate}
\end{thm}

\begin{df}
Let $B_1$ and $B_2$ be crystals. A {\it strict morphism} of crystals
$\psi:B_1\lar B_2$ is a map
$\psi:B_1\sqcup\{0\} \lar B_2\sqcup\{0\}$
satisfying: 
$\psi(0)=0$, $\psi(B_1)\subset B_2$,
$\psi$ commutes with all $\eit$ and $\fit$
and
\[
\hspace{-30pt}\wt(\psi(b))=\wt(b),\q \vep_i(\psi(b))=\vep_i(b),\q
  \vp_i(\psi(b))=\vp_i(b)
\text{ for any }b\in B_1.
\]
In particular, 
a bijective strict morphism is called an 
{\it isomorphism of crystals}. 
\end{df}

\begin{ex}
\label{ex-tlm}
If $(L,B)$ is a crystal base, then $B$ is a crystal.
Hence, for the crystal base $(L(\ify),B(\ify))$
of the nilpotent subalgebra $\uqm$ of 
the quantum algebra $\uq$, 
$B(\ify)$ is a crystal. 
\end{ex}
\begin{ex}
\label{tlm}
For $\lm\in P$, set $T_\lm:=\{t_\lm\}$. We define a crystal
structure on $T_\lm$ by 
\[
 \eit(t_\lm)=\fit(t_\lm)=0,\q\vep_i(t_\lm)=
\vp_i(t_\lm)=-\ify,\q \wt(t_\lm)=\lm.
\]
\end{ex}
\begin{df}
For a crystal $B$, a colored oriented graph
structure is associated with $B$ by 
\[
 b_1\mapright{i}b_2\Longleftrightarrow 
\fit b_1=b_2.
\]
We call this graph a {\it crystal graph}
of $B$.
\end{df}

\subsection{Affine weights}
\label{aff-wt}

Let $\ge$ be an affine Lie algebra. 
The sets $\mathfrak t$, 
$\{\al_i\}_{i\in I}$ 
and $\{\al^\vee_i\}_{i\in I}$ be as in \ref{KM}. 
We take ${\rm dim}\mathfrak t=\sharp I+1$.
Let $\del\in Q_+$ be the unique element 
satisfying $\{\lm\in Q|\lan \al^\vee_i,\lm\ran=0
\text{ for any }i\in I\}=\bbZ\del$
and ${\bf c}\in \ge$ be the canonical central element
satisfying $\{h\in Q^\vee|\lan h,\al_i\ran=0
\text{ for any }i\in I\}=\bbZ c$.
We write (\cite[6.1]{Kac})
\[
{\bf c}=\sum_i a_i^\vee \al^\vee_i,\qq
\del=\sum_i a_i\al_i.
\]
Let $(\q,\q)$ be the non-degenerate
$W$-invariant symmetric bilinear form on $\mathfrak t^*$
normalized by $(\del,\lm)=\lan {\bf c},\lm\ran$
for $\lm\in\frak t^*$.
Let us set $\tt^*_{\rm cl}:=\tt^*/\bbC\del$ and let
${\rm cl}:\tt^*\longrightarrow \tt^*_{\rm cl}$
be the canonical projection. 
Here we have 
$\tt^*_{\rm cl}\cong \oplus_i(\bbC \al^\vee_i)^*$.
Set $\tt^*_0:=\{\lm\in\tt^*|\lan {\bf c},\lm\ran=0\}$,
$(\tt^*_{\rm cl})_0:={\rm cl}(\tt^*_0)$. 
Since $(\del,\del)=0$, we have a positive-definite
symmetric form on $\tt^*_{\rm cl}$ 
induced by the one on 
$\tt^*$. 
Let $\Lm_i\in \tt^*_{\rm cl}$ $(i\in I)$ be a classical 
weight such that $\lan \al^\vee_i,\Lm_j\ran=\del_{i,j}$, which 
is called a fundamental weight.
We choose 
$P$ so that $P_{\rm cl}:={\rm cl}(P)$ coincides with 
$\oplus_{i\in I}\bbZ\Lm_i$ and 
we call $P_{\rm cl}$ a 
{\it classical weight lattice}.

\subsection{Definitions of perfect crystal and its limit}
\label{def-perfect}

Let $\ge$ be an affine Lie algebra, $P_{cl}$ be
a classical weight lattice as above and set 
$(P_{cl})^+_l:=\{\lm\in P_{cl}|
\lan c,\lm\ran=l,\,\,\lan \al^\vee_i,\lm\ran\geq0\}$ 
$(l\in\ZZ_{>0})$.
\begin{df}
\label{perfect-def}
A crystal $B$ is a {\it perfect} of level $l$ if 
\begin{enumerate}
\item
$B\ot B$ is connected as a crystal graph.
\item
There exists $\lm_0\in P_{\rm cl}$ such that 
\[
 \wt(B)\subset \lm_0+\sum_{i\ne0}\ZZ_{\leq0}
{\rm cl}(\al_i),\qq
\sharp B_{\lm_0}=1
\]
\item There exists a finite-dimensional 
$U'_q(\ge)$-module $V$ with a
crystal pseudo-base $B_{ps}$ 
such that $B\cong B_{ps}/{\pm1}$
\item
The maps 
$\vep,\vp:B^{min}:=\{b\in B|\lan c,\vep(b)\ran=l\}
\mapright{}(P_{\rm cl}^+)_l$ are bijective, where 
$\vep(b):=\sum_i\vep_i(b)\Lm_i$ and 
$\vp(b):=\sum_i\vp_i(b)\Lm_i$.
\end{enumerate}
\end{df}

Let $\{B_l\}_{l\geq1}$ be a family of 
perfect crystals of level $l$ and set 
$J:=\{(l,b)|l>0,\,b\in B^{min}_l\}$.
\begin{df}
\label{def-limit}
A crystal $B_\ify$ with an element $b_\ify$ is called a
{\it limit of $\{B_l\}_{l\geq1}$}
if 
\begin{enumerate}
\item
$\wt(b_\ify)=\vep(b_\ify)=\vp(b_\ify)=0$.
\item
For any $(l,b)\in J$, there exists an
embedding of 
crystals:
\begin{eqnarray*}
 f_{(l,b)}:&
T_{\vep(b)}\ot B_l\ot T_{-\vp(b)}\hookrightarrow
B_\ify\\
&t_{\vep(b)}\ot b\ot t_{-\vp(b)}\mapsto b_\ify
\end{eqnarray*}
\item
$B_\ify=\bigcup_{(l,b)\in J} {\rm Im}f_{(l,b)}$.
\end{enumerate}
\end{df}
\noindent
As for the crystal $T_\lm$, see Example \ref{tlm}.
If a limit exists for a family $\{B_l\}$, 
we say that $\{B_l\}$
is a {\it coherent family} of perfect crystals.

The following is one of the most 
important properties of limit of perfect crystals.
\begin{pro}
Let $B(\ify)$ be the crystal as in 
Example \ref{ex-tlm}. Then we have
the following isomorphism of crystals:
\[
B(\ify)\ot B_\ify\mapright{\sim}B(\ify).
\]
\end{pro}

\section{Perfect Crystals of type $\TY(D,3,4)$}
\label{perf}

In this section, we review the family of 
perfect crystals of type 
$\TY(D,3,4)$ and its limit(\cite{KMOY}).

We fix the data for $D_4^{(3)}$.
Let 
$\{\alpha_0, \alpha_1, \alpha_2\}$, 
$\{\al^\vee_0, \al^\vee_1, \al^\vee_2\}$ and 
$\{\Lm_0, \Lm_1, \Lm_2\}$ be the set of 
simple roots, simple coroots and fundamental weights, respectively.
The Cartan matrix $A=(a_{ij})_{i,j=0,1,2}$
is given by
\[A=
\left(
\begin{array}{rrr}
2  & -1 & 0  \\
-1 & 2  & -3 \\
0  & -1 & 2
\end{array}
\right),
\]
and its Dynkin diagram is as follows.
\[\SelectTips{cm}{}
\xymatrix{
*{\bigcirc}<3pt> \ar@{-}[r]_<{0} 
& *{\bigcirc}<3pt> \ar@3{<-}[r]_<{1}
& *{\bigcirc}<6pt>\ar@{}_<{\,\,\,\,\,\,2}
}
\]
The standard null root $\delta$ 
and the canonical central element $c$ are 
given by
\[
\delta=\alpha_0+2\alpha_1+\alpha_2
\quad\text{and}\quad c=\al^\vee_0+2\al^\vee_1+3\al^\vee_2, 
\]
where 
$\al_0=2\Lm_0-\Lm_1+\del,\q
\al_1=-\Lm_0+2\Lm_1-\Lm_2,\q
\al_2=-3\Lm_1+2\Lm_2.$

For a positive integer $l$ we introduce 
$\TY(D,3,4)$-crystals $B_l$ and $B_\ify$ as 
\begin{eqnarray*}
&&B_l=\left\{
b=(b_1,b_2,b_3,{\bar b}_3,{\bar b}_2,{\bar b}_1)
\in(\Zn)^6
\left\vert
\begin{array}{l}
b_3\equiv\bar{b}_3\;(\text{mod }2), \\
\sum_{i=1,2} (b_i+{\bar b}_i)+\frac{b_3+{\bar b}_3}{2}
\leq l
\end{array}
\right.
\right\},\\
&&B_\ify=\left\{
b=(b_1,b_2,b_3,{\bar b}_3,{\bar b}_2,{\bar b}_1)
\in(\ZZ)^6
\left\vert
\begin{array}{l}
b_3\equiv\bar{b}_3\;(\text{mod }2), \\
\sum_{i=1,2} (b_i+{\bar b}_i)+\frac{b_3+{\bar b}_3}{2}
\in\ZZ
\end{array}
\right.
\right\}. 
\end{eqnarray*}
Now we describe the explicit crystal structures
of $B_l$ and $B_\ify$. 
Indeed, most of them coincide with 
each other except
for $\vep_0$ and $\vp_0$.
In the rest of this section, we use the following 
convention: 
$(x)_+=\max(x,0)$.
\begin{align*}
{\tilde e}_1 b=&
\begin{cases}
(\ldots,{\bar b}_2 +1,{\bar b}_1 -1) 
& \text{if ${\bar b}_2 -{\bar b}_3 \geq (b_2 -b_3)_+$}, 
\\  
(\ldots,b_3 +1,{\bar b}_3 -1,\ldots) 
& \text{if ${\bar b}_2 -{\bar b}_3 <0\leq b_3 -b_2$}, 
\\ 
(b_1 +1,b_2 -1,\ldots) 
& \text{if $({\bar b}_2 -{\bar b}_3)_+ <b_2 -b_3$},
\end{cases}
\\
{\tilde f}_1 b=&
\begin{cases}
(b_1 -1,b_2 +1,\ldots) 
& \text{if $({\bar b}_2 -{\bar b}_3)_+ \leq b_2 -b_3$}, 
\\
(\ldots,b_3 -1,{\bar b}_3 +1,\ldots) 
& \text{if ${\bar b}_2 -{\bar b}_3 \leq 0<b_3 -b_2$}, 
\\
(\ldots,{\bar b}_2 -1,{\bar b}_1 +1) 
& \text {if ${\bar b}_2 -{\bar b}_3 >(b_2 -b_3)_+$},
\end{cases}
\\
{\tilde e}_2 b=&
\begin{cases}
(\ldots,{\bar b}_3 +2,{\bar b}_2 -1,\ldots) 
& \text{if ${\bar b}_3 \geq b_3$}, 
\\
(\ldots,b_2 +1,b_3 -2,\ldots) 
& \text{if ${\bar b}_3 <b_3$},
\end{cases}
\\
{\tilde f}_2 b=&
\begin{cases}
(\ldots,b_2 -1,b_3 +2,\ldots) 
& \text{if ${\bar b}_3 \leq b_3$}, 
\\
(\ldots,{\bar b}_3 -2,{\bar b}_2 +1,\ldots) 
& \text{if ${\bar b}_3 >b_3$},
\end{cases}
\end{align*}
\begin{align*}
\vep_1(b)=&{\bar b}_1+({\bar b}_3-{\bar b}_2+(b_2-b_3)_+)_+,
\qq\varphi_1(b)=
b_1+(b_3-b_2+({\bar b}_2-{\bar b}_3)_+)_+,
\\
\vep_2(b)=&{\bar b}_2+\frac{1}{2}(b_3-{\bar b}_3)_+,
\qq\varphi_2(b)=b_2+\frac{1}{2}({\bar b}_3-b_3)_+,\\
\vep_0(b)=&
\begin{cases}
l-s(b)+\max \,A-(2z_1+z_2+z_3+3z_4)&b\in B_l,\\
-s(b)+\max \,A-(2z_1+z_2+z_3+3z_4)&b\in B\ify.
\end{cases}\\
\vp_0(b)=&
\begin{cases}
l-s(b)+\max\, A&b\in B_l,\\
-s(b)+\max\, A&b\in B_\ify,
\end{cases}
\end{align*}
where
\begin{equation} \label{def s(b)}
s(b)=b_1+b_2+\frac{b_3+{\bar b}_3}{2}+{\bar b}_2+{\bar b}_1. 
\end{equation}
\begin{equation}
 \label{z1-4}
z_1={\bar b}_1-b_1, \quad 
z_2={\bar b}_2 -{\bar b}_3, \quad 
z_3=b_3-b_2, \quad 
z_4=({\bar b}_3-b_3)/2,
\end{equation}
\begin{equation} \label{A}
A=(0,z_1,z_1+z_2,z_1+z_2+3z_4,z_1+z_2+z_3+3z_4,2z_1+z_2+z_3+3z_4)
\end{equation}
For $b\in B_l$ 
if $\et{i}b$ or $\ft{i}b$ does not belong to 
$B_l$, namely, if $b_j$ or $\bar{b}_j$
for some $j$ becomes negative, 
we understand it to be $0$.

Let us see the actions of $\et{0}$ and $\ft{0}$.
We shall
consider the 
conditions ($E_1$)-($E_6$) and ($F_1$)-($F_6$) 
(\cite{KMOY}). 
\begin{eqnarray*}
&&(E_1)\quad
z_1+z_2+z_3+3z_4<0, z_1+z_2+3z_4<0, 
z_1+z_2<0, z_1<0,
\\
&&(E_2)\quad
z_1+z_2+z_3+3z_4<0, z_2+3z_4<0, z_2<0, 
z_1\geq 0,
\\
&&(E_3)\quad 
z_1+z_3+3z_4<0, z_3+3z_4<0, z_4<0, z_2\geq 0, 
z_1+z_2\geq 0,
\\
&&(E_4)\quad
z_1+z_2+3z_4\geq 0, z_2+3z_4\geq 0, z_4\geq 0, 
z_3<0, z_1+z_3<0,
\\
&&(E_5)\quad
z_1+z_2+z_3+3z_4\geq 0, z_3+3z_4\geq 0, 
z_3\geq 0, z_1<0, 
\\
&&(E_6)\quad
z_1+z_2+z_3+3z_4\geq 0, z_1+z_3+3z_4\geq 0, 
z_1+z_3\geq 0, z_1\geq 0.
\end{eqnarray*}
($F_i$) ($1\le i\le 6$) is obtained 
from ($E_i$) by 
replacing $\geq$ (resp. $<$) with 
$>$ (resp. $\leq$). We define
\begin{align*}
{\tilde e}_0 b=&
\begin{cases}
\mathscr{E}_1 b:=
(b_1 -1,\ldots) 
& \text{if ($E_1$)}, 
\\
\mathscr{E}_2 b:=
(\ldots,b_3 -1,{\bar b}_3 -1,\ldots,{\bar b}_1 +1) 
& \text{if ($E_2$)}, 
\\
\mathscr{E}_3 b:=
(\ldots,b_3 -2,\ldots,{\bar b}_2 +1,\ldots) 
& \text{if ($E_3$)}, 
\\ 
\mathscr{E}_4 b:=
(\ldots,b_2 -1,\ldots,{\bar b}_3 +2,\ldots) 
& \text{if ($E_4$)}, 
\\
\mathscr{E}_5 b:=
(b_1 -1,\ldots,b_3 +1,{\bar b}_3 +1,\ldots) 
& \text{if ($E_5$)}, 
\\
\mathscr{E}_6 b:=
(\ldots,{\bar b}_1 +1) 
& \text{if ($E_6$)},
\end{cases}
\\
{\tilde f}_0 b=&
\begin{cases}
\mathscr{F}_1 b:=(b_1 +1,\ldots) 
& \text{if ($F_1$)}, 
\\ 
\mathscr{F}_2 b:=
(\ldots,b_3 +1,{\bar b}_3 +1,\ldots,{\bar b}_1 -1) 
& \text{if ($F_2$)}, 
\\
\mathscr{F}_3 b:=
(\ldots,b_3 +2,\ldots,{\bar b}_2 -1,\ldots) 
& \text{if ($F_3$)}, 
\\
\mathscr{F}_4 b:=
(\ldots,b_2 +1,\ldots,{\bar b}_3 -2,\ldots) 
& \text{if ($F_4$)}, 
\\
\mathscr{F}_5 b:=
(b_1 +1,\ldots,b_3 -1,{\bar b}_3 -1,\ldots) 
& \text{if ($F_5$)}, 
\\
\mathscr{F}_6 b:=(\ldots,{\bar b}_1 -1)
& \text{if ($F_6$)}.
\end{cases}
\end{align*}
The following is one of the main results in 
\cite{KMOY}:
\begin{thm}[\cite{KMOY}]
\begin{enumerate}
\item
The $\TY(D,3,4)$-crystal $B_l$ 
is a perfect crystal of level $l$. 
\item
The family of the perfect crystals 
$\{B_l\}_{\l\geq1}$ forms a 
coherent family and the crystal $B_\ify$ 
is its limit with the vector 
$b_\ify=(0,0,0,0,0,0)$.
\end{enumerate}
\end{thm}
As was shown in \cite{KMOY}, 
the minimal elements are given
\[
(B_l)_{\min}=
\{(\alpha,\beta,\beta,\beta,\beta,\alpha)\,|\,
\al,\beta\in\ZZ_{\geq0}, 2\al+3\beta\leq l\}.
\] 
Let $J=\{(l,b)\,|\,l\in\ZZ_{\geq1},
b\in(B_l)_{\min}\}$ and the maps $\vep,\,\vp
:(B_l)_{\min}
\to (P^+_{\rm cl})_l$ be as in Sect.3.
Then we have 
$\wt b_\infty=0$ and $\vep_i(b_\infty)
=\vp_i(b_\infty)=0$ for $i=0,1,2$.

For $(l,b_0)\in J$, since $\vep(b_0)=\vp(b_0)$, 
one can set $\lm=\vep(b_0)=\vp(b_0)$. For 
$b=(b_1,b_2,b_3,\bar{b}_3,\bar{b}_2,\bar{b}_1)\in B_l$ we define a map
\[
f_{(l,b_0)}:\;T_{\lm}\ot B_l\ot 
B_{-\lm}\longrightarrow B_\infty
\]
by
\[
f_{(l,b_0)}(t_{\lm}\ot b\ot t_{-\lm})=b'=(\nu_1,\nu_2,\nu_3,\bar{\nu}_3,\bar{\nu}_2,\bar{\nu}_1)
\]
where
$b_0=(\al,\beta,\beta,\beta,\beta,\al)$, and
\begin{align*}
\nu_1&=b_1-\alpha, & \bar{\nu}_1&=\bar{b}_1-\alpha, \\
\nu_j&=b_j-\beta, & \bar{\nu}_j&=\bar{b}_j-\beta\;(j=2,3).
\end{align*}
Finally, we obtain 
$B_\infty=\bigcup_{(l,b)\in J}
\mbox{Im}\,f_{(l,b)}$

\renewcommand{\thesection}{\arabic{section}}
\section{Fundamental Representation for 
$\TY(G,1,2)$}
\setcounter{equation}{0}
\renewcommand{\theequation}{\thesection.\arabic{equation}}

\subsection{Fundamental representation 
$W(\varpi_1)$}
\label{fundamental}

Let $c=\sum_{i}a_i^\vee \al^\vee_i$ be the canonical
central element in an affine Lie algebra $\ge$
(see \cite[6.1]{Kac}), 
$\{\Lm_i|i\in I\}$ the set of fundamental 
weight as in the previous section
and $\varpi_1:=\Lm_1-a^\vee_1\Lm_0$ the
(level 0)fundamental weight.
Let $W(\varpi_1)$ be the fundamental representation 
of $\uqp$
associated with $\varpi_1$ (\cite{K0}).

By \cite[Theorem 5.17]{K0}, $W(\varpi_1)$ is a
finite-dimensional irreducible integrable 
$\uqp$-module and has a global basis
with a simple crystal. Thus, we can consider 
the specialization $q=1$ and obtain the 
finite-dimensional $\ge$-module $W(\varpi_1)$, 
which we call a fundamental representation
of $\ge$ and use the same notation as above.

We shall present the explicit form of 
$W(\varpi_1)$ for $\ge=\TY(G,1,2)$.
\subsection{$W(\varpi_1)$ for $\TY(G,1,2)$}
The Cartan matrix $A=(a_{i,j})_{i,j=0,1,2}$ of type 
$\TY(G,1,2)$ is:
\[
 A=\begin{pmatrix}2&-1&0\\
-1&2&-1\\0&-3&2
\end{pmatrix}.
\]
Then the simple roots are 
\[
 \al_0=2\Lm_0-\Lm_1+\del,\q
\al_1=-\Lm_0+2\Lm_1-3\Lm_2,\q
\al_2=-\Lm_1+2\Lm_2, 
\]
and the Dynkin diagram is:
\[\SelectTips{cm}{}
\xymatrix{
*{\bigcirc}<3pt> \ar@{-}[r]_<{0} 
& *{\bigcirc}<3pt> \ar@3{->}[r]_<{1}
& *{\bigcirc}<6pt>\ar@{}_<{\,\,\,\,\,\,2}
}
\]

The $\ge$-module $W(\varpi_1)$ is a 15 dimensional 
module with the basis,
\[
 \{\fsquare(5mm,i),\fsquare(5mm,\ovl i),
\emptyset, \fsquare(5mm,0_1),
\fsquare(5mm,0_2) \,\,\vert \,\,i=1,\cd,6\}.
\]
The following description of $W(\varpi_1)$ 
slightly differs from \cite{Y}.
\begin{eqnarray*}
&&\hspace{-30pt}{\rm wt}(\fsquare(5mm,1))=\Lm_1-2\Lm_0,\,\,
{\rm wt}(\fsquare(5mm,2))=-\Lm_0-\Lm_1+3\Lm_2,\,\,
{\rm wt}(\fsquare(5mm,3))=-\Lm_0+\Lm_2,\\
&&\hspace{-30pt}{\rm wt}(\fsquare(5mm,4))=-\Lm_0+\Lm_1-\Lm_2,\,\,
{\rm wt}(\fsquare(5mm,5))=-\Lm_1+2\Lm_2,\,\,
{\rm wt}(\fsquare(5mm,6))=-\Lm_0+2\Lm_1-3\Lm_2,\\
&&\hspace{-30pt}{\rm wt}(\fsquare(5mm,\ovl i))=
-{\rm wt}(\fsquare(5mm,i))\,\,(i=1,\cd,6),\,\,
{\rm wt}(\fsquare(5mm,0_1))=
{\rm wt}(\fsquare(5mm,0_2))=
{\rm wt}(\emptyset)=0.
\end{eqnarray*}
\def\bv#1{\fsquare(5mm,#1)}
The actions of $e_i$ and $f_i$ on these basis vectors
are given as follows:
\begin{eqnarray*}
&&\hspace{-30pt}
f_0\left(\bv{0_2},\bv{\ovl 6},\bv{\ovl 4},\bv{\ovl 3},\bv{\ovl 2},
\bv{\ovl 1},\emptyset\right)
=\left(\bv1,\bv2,\bv3,\bv4,\bv6,\emptyset,2\bv1\right),\\
&&\hspace{-30pt}e_0\left(\bv1,\bv2,\bv3,\bv4,\bv6,\bv{0_2},\emptyset\right)
=\left(\emptyset,\bv{\ovl 6},\bv{\ovl 4},\bv{\ovl 3},
\bv{\ovl 2},\bv{\ovl 1},2\bv{\ovl 1}\right),\\
&&\hspace{-30pt}f_1\left(\bv1,\bv4,\bv6,\bv{0_1},\bv{0_2},
\bv{\ovl 5},\bv{\ovl 2},\emptyset\right)=
\left(\bv2,\bv5,\bv{0_2},3\bv{\ovl 6},2\bv{\ovl 6},
\bv{\ovl 4},\bv{\ovl 1},\bv{\ovl 6}\right),\\
&&\hspace{-30pt}e_1\left(\bv2,\bv5,\bv{0_1},\bv{0_2},\bv{\ovl 6},
\bv{\ovl 4},\bv{\ovl 1},\emptyset\right)=
\left(\bv1,\bv4,3\bv6,2\bv6,\bv{0_2},\bv{\ovl 5},
\bv{\ovl 2},\bv6\right),\\
&&\hspace{-30pt}f_2\left(\bv2,\bv3,\bv4,\bv5,\bv{0_1},\bv{0_2},
\bv{\ovl 6},\bv{\ovl 4},\bv{\ovl 3}\right)\\
&&\qq\qq\qq\qq\qq\qq=
\left(\bv3,2\bv4,3\bv6,\bv{0_1},2\bv{\ovl 5},\bv{\ovl 5},
\bv{\ovl 4},2\bv{\ovl 3},3\bv{\ovl 2}\right),\\
&&\hspace{-30pt}e_2\left(\bv3,\bv4,\bv6,\bv{0_1},\bv{0_2},\bv{\ovl 5},
\bv{\ovl 4},\bv{\ovl 3},\bv{\ovl 2}\right)\\
&&\qq\qq\qq\qq\qq\qq=
\left(3\bv2,2\bv3,\bv4,2\bv5,\bv5,\bv{0_1},
3\bv{\ovl 6},2\bv{\ovl 4},\bv{\ovl 3}\right),
\end{eqnarray*}
where we give non-trivial actions only.

\renewcommand{\thesection}{\arabic{section}}
\section{Affine Geometric Crystal $\cV_1(\TY(G,1,2))$}
\setcounter{equation}{0}
\renewcommand{\theequation}{\thesection.\arabic{equation}}

Let us review the construction of the 
affine geometric crystal $\cV(\TY(G,1,2))$ in $W(\varpi_1)$ following \cite{N3}.

For $\xi\in (\frt^*_{\rm cl})_0$, let $t(\xi)$ be the 
shift as in \cite[Sect 4]{K0}.
Then we have 
\begin{eqnarray*}
&& t(\wtil\varpi_1)=s_0s_1s_2s_1s_2s_1=:w_1,\\
&& t(\text{wt}(\bv{\ovl 2}))=s_2s_1s_2s_1s_0s_1=:w_2,
\end{eqnarray*}
Associated with these Weyl group elements $w_1$ and $w_2$,
we define algebraic varieties $\cV_1=\cV_1(\TY(G,1,2))$ and 
$\cV_2=\cV_2(\TY(G,1,2))\subset W(\varpi_1)$ respectively:
\begin{eqnarray*}
&&\hspace{-30pt}\cV_1:=\{v_1(x)
:=Y_0(x_0)Y_1(x_1)Y_2(x_2)Y_1(x_3)Y_2(x_4)Y_1(x_5)
\bv1\,\,\vert\,\,x_i\in\bbC^\times,(0\leq i\leq 5)\},\\
&&\hspace{-30pt}\cV_2:=\{v_2(y):=
Y_2(y_2)Y_1(y_1)Y_2(y_4)Y_1(y_3)Y_0(y_0)Y_1(y_5)
\bv{\ovl 2}\,\,\vert\,\,y_i\in\bbC^\times,
(0\leq i\leq 5)\}.
\end{eqnarray*}
Owing to the explicit forms of $f_i$'s on $W(\varpi_1)$
as above, we have $f_0^3=0$, $f_1^3=0$ and $f_2^4=0$ 
and then 
\[
Y_i(c)=(1+\frac{f_i}{c}+\frac{f_i^2}{2c^2})\al_i^\vee(c)
\,\,(i=0,1),\q
Y_2(c)=(1+\frac{f_2}{c}+\frac{f_2^2}{2c^2}
+\frac{f_2^3}{6c^3})\al_2^\vee(c).
\]
We get explicit forms of $v_1(x)\in\cV_1$ 
and $v_2(y)\in\cV_2$ as in \cite{N3}:
\begin{eqnarray*}
&&v_1(x)=\sum_{1\leq i\leq 6}\left(X_i\bv{i}+X_{\ovl i}
\bv{\ovl i}\right)+X_{0_1}\bv{0_1}+X_{0_2}\bv{0_2}+X_\emptyset
\emptyset,\\
&&v_2(y)=\sum_{1\leq i\leq 6}\left(Y_i\bv{i}+Y_{\ovl i}
\bv{\ovl i}\right)+Y_{0_1}\bv{0_1}+Y_{0_2}\bv{0_2}+Y_\emptyset
\emptyset.
\end{eqnarray*}
where the rational functions $X_i$'s and $Y_i$'s are all
positive (as for their explicit forms, see \cite{N3})
and then 
we get the positive birational isomorphism 
$\ovl\sigma:\cV_1\longrightarrow \cV_2$ ($v_1(x)\mapsto v_2(y)$)
and its inverse $\ovl\sigma^{-1}$ is also positive.
The actions of $e_0^c$ on $v_2(y)$ 
(respectively $\gamma_0(v_2(y))$ and 
$\vep_0(v_2(y)))$ are 
induced from the ones on 
$Y_2(y_2)Y_1(y_1)Y_2(y_4)Y_1(y_3)Y_0(y_0)Y_1(y_5)$ 
as an element of the geometric crystal $\cV_2$.
We define the action $e_0^c$ on $v_1(x)$ by
\begin{equation}
e_0^cv_1(x)=\ovl\sigma^{-1}\circ e_0^c\circ
\ovl\sigma(v_1(x))).
\label{e0}
\end{equation}
We also define $\gamma_0(v_1(x))$ 
and $\vep_0(v_1(x))$ by 
\begin{equation}
\gamma_0(v_1(x))=\gamma_0(\ovl\sigma(v_1(x))),\qq
\vep_0(v_1(x)):=\vep_0(\ovl\sigma(v_1(x))).
\label{wt0}
\end{equation}
\begin{thm}[\cite{N3}]
Together with (\ref{e0}), 
(\ref{wt0}) on $\cV_1$, we obtain a
positive affine geometric crystal $\chi:=
(\cV_1,\{e_i\}_{i\in I},
\{\gamma_i\}_{i\in I},\{\vep_i\}_{i\in I})$
$(I=\{0,1,2\})$, whose explicit form is as follows:
first we have $e_i^c$, $\gamma_i$ and $\vep_i$
for $i=1,2$ from the formula (\ref{eici}), 
(\ref{vep-i})
and (\ref{gamma-i}).
\begin{eqnarray*}
&&\hspace{-30pt}
e_1^c(v_1(x))=v_1(x_0,\cC_1x_1,x_2,
\cC_3x_3,x_4,\cC_5x_5),\,
e_2^c(v_1(x))=v_1(x_0,x_1,
\cC_2x_2,x_3,\cC_4x_4,x_5),\\
&&\text{where}\\
&&\cC_1=\frac{\frac{c\,{x_0}}{{x_1}} 
+ \frac{{x_0}\,{{x_2}}^3}{{{x_1}}^2\,{x_3}} 
+\frac{{x_0}\,{{x_2}}^3\,{{x_4}}^3}{{{x_1}}^2\,
{{x_3}}^2\,{x_5}}}{\frac{{x_0}}{{x_1}} 
+ \frac{{x_0}\,{{x_2}}^3}{{{x_1}}^2\,{x_3}} + 
\frac{{x_0}\,{{x_2}}^3\,
{{x_4}}^3}{{{x_1}}^2\,{{x_3}}^2\,{x_5}}},\q
\cC_3=\frac{\frac{c\,{x_0}}{{x_1}} 
+ \frac{c\,{x_0}\,{{x_2}}^3}{{{x_1}}^2\,{x_3}} + 
\frac{{x_0}\,{{x_2}}^3\,{{x_4}}^3}{{{x_1}}^2
\,{{x_3}}^2\,{x_5}}}{
\frac{c\,{x_0}}{{x_1}} 
+ \frac{{x_0}\,{{x_2}}^3}{{{x_1}}^2\,{x_3}} + 
\frac{{x_0}\,{{x_2}}^3\,{{x_4}}^3}
{{{x_1}}^2\,{{x_3}}^2\,{x_5}}},\\
&&\cC_5=\frac{c\,\left( \frac{{x_0}}{{x_1}} + 
      \frac{{x_0}\,{{x_2}}^3}{{{x_1}}^2\,{x_3}} + 
\frac{{x_0}\,{{x_2}}^3\,{{x_4}}^3}{{{x_1}}^2
\,{{x_3}}^2\,{x_5}} \right)}{\frac{c\,{x_0}}{{x_1}} 
+ \frac{c\,{x_0}\,{{x_2}}^3}{{{x_1}}^2\,{x_3}} + 
\frac{{x_0}\,{{x_2}}^3\,{{x_4}}^3}
{{{x_1}}^2\,{{x_3}}^2\,{x_5}}},\,
\cC_2=\frac{\frac{c\,{x_1}}{{x_2}} 
+ \frac{{x_1}\,{x_3}}{{{x_2}}^2\,{x_4}}}
{\frac{{x_1}}{{x_2}} + \frac{{x_1}\,{x_3}}{{{x_2}}^2\,{x_4}}},\,
\cC_4=\frac{c\,\left( \frac{{x_1}}{{x_2}} + 
\frac{{x_1}\,{x_3}}{{{x_2}}^2\,{x_4}} \right) }{\frac{c\,{x_1}}
{{x_2}} + \frac{{x_1}\,{x_3}}{{{x_2}}^2\,{x_4}}},\\
&&\vep_1(v_1(x))={\frac{{x_0}}{{x_1}} 
+ \frac{{x_0}\,{{x_2}}^3}{{{x_1}}^2\,{x_3}} + 
\frac{{x_0}\,{{x_2}}^3\,
{{x_4}}^3}{{{x_1}}^2\,{{x_3}}^2\,{x_5}}},\q
\vep_2(v_1(x))={\frac{{x_1}}{{x_2}} 
+ \frac{{x_1}\,{x_3}}{{{x_2}}^2\,{x_4}}},\\
&&\gamma_1(v_1(x))=\frac{x_1^2x_3^2x_5^2}{x_0x_2^3x_4^3},\q
\gamma_2(v_1(x))=\frac{x_2^2x_4^2}{x_1x_3x_5}.
\end{eqnarray*}
We also have $e_0^c$, $\vep_0$ and $\gamma_0$ on $v_1(x)$
as:
\begin{eqnarray*}
&&e_0^c(v_1(x))=v_1(\frac{D}{c\cdot E}x_0,\frac{F}{c\cdot E}x_1,
\frac{G}{c\cdot E}x_2,\frac{D\cdot H}{c^2\cdot E\cdot F}x_3,
\frac{D}{c\cdot G}x_4,\frac{D}{c\cdot H}x_5),\\
&&\vep_0(v_1(x))=\frac{E}{{{x_0}}^3\,{{x_2}}^3\,{x_3}},\qq
\gamma_0(v_1(x))=\frac{x_0^2}{x_1x_3x_5},\\
&&\text{where}\\
&&D=c^2\,{{x_0}}^2\,{{x_2}}^3\,{x_3} + 
  {x_1}\,{{x_2}}^3\,{{x_3}}^2\,{x_5} + 
  c\,{x_0}\,( {x_1}\,{{x_3}}^3 + 
     3\,{x_1}\,{x_2}\,{{x_3}}^2\,{x_4}\\
&&\qq\qq\qq + 3\,{x_1}\,{{x_2}}^2\,{x_3}\,{{x_4}}^2
+ {{x_2}}^3\,\left( {{x_3}}^2 + {x_1}\,{{x_4}}^3 + 
        {x_1}\,{x_3}\,{x_5} \right) ),\\
&&E={{x_0}}^2\,{{x_2}}^3\,{x_3} 
+ {x_1}\,{{x_2}}^3\,{{x_3}}^2\,{x_5} + 
  {x_0}\,\left( {x_1}\,{{x_3}}^3 + 
     3\,{x_1}\,{x_2}\,{{x_3}}^2\,{x_4} + 
     3\,{x_1}\,{{x_2}}^2\,{x_3}\,{{x_4}}^2 \right.\\
&&\qq \left.+  {{x_2}}^3\,\left( {{x_3}}^2 + {x_1}\,{{x_4}}^3 + 
        {x_1}\,{x_3}\,{x_5} \right)  \right),\\
&&F=c\,{{x_0}}^2\,{{x_2}}^3\,{x_3} + 
  {x_1}\,{{x_2}}^3\,{{x_3}}^2\,{x_5} + 
  {x_0}\,( c\,{x_1}\,{{x_3}}^3 + 
     3\,c\,{x_1}\,{x_2}\,{{x_3}}^2\,{x_4}\\
&&\qq\qq\qq + 3\,c\,{x_1}\,{{x_2}}^2\,{x_3}\,{{x_4}}^2
+ {{x_2}}^3\,\left( {{x_3}}^2 + c\,{x_1}\,{{x_4}}^3 + 
        c\,{x_1}\,{x_3}\,{x_5} \right)),\\
&&G=c\,{{x_0}}^2\,{{x_2}}^3\,{x_3} + 
  {x_1}\,{{x_2}}^3\,{{x_3}}^2\,{x_5} + 
  {x_0}\,( {x_1}\,{{x_3}}^3 + 
     \left( 2 + c \right) \,{x_1}\,{x_2}\,{{x_3}}^2\,{x_4}\\
&&\qq\qq\qq + \left( 1 + 2\,c \right) 
\,{x_1}\,{{x_2}}^2\,{x_3}\,{{x_4}}^2
+ {{x_2}}^3\,\left( {{x_3}}^2 + c\,{x_1}\,{{x_4}}^3 + 
        c\,{x_1}\,{x_3}\,{x_5}\right)),\\
&&H=c\,{{x_0}}^2\,{{x_2}}^3\,{x_3} + 
  {x_1}\,{{x_2}}^3\,{{x_3}}^2\,{x_5} + 
  {x_0}\,( {x_1}\,{{x_3}}^3 + 
     3\,{x_1}\,{x_2}\,{{x_3}}^2\,{x_4} \\
&&\qq\qq\qq
+ 3\,{x_1}\,{{x_2}}^2\,{x_3}\,{{x_4}}^2
+ {{x_2}}^3\,\left( {{x_3}}^2 + {x_1}\,{{x_4}}^3 + 
        c\,{x_1}\,{x_3}\,{x_5} \right)).
\end{eqnarray*}
\end{thm}

\renewcommand{\thesection}{\arabic{section}}
\section{Ultra-discretization}
\setcounter{equation}{0}
\renewcommand{\theequation}{\thesection.\arabic{equation}}

We denote the positive structure on $\chi$ as in the 
previous section by 
$\theta:T'\seteq(\bbC^\times)^6 \longrightarrow \cV_1$. 
Then by Corollary \ref{cor-posi}
we obtain the ultra-discretization 
${\mathcal UD}(\chi,T',\theta)$, 
which is a Kashiwara's crystal. 
Now we show that the conjecture in \cite{N3} 
is correct and it turns 
out to be the following theorem.
\begin{thm}
\label{ultra-d}
The crystal ${\mathcal UD}(\chi,T',\theta)$ as above 
is isomorphic to the crystal $B_\ify$ of 
type $\TY(D,3,4)$ as in Sect.\ref{perf}.
\end{thm}
In order to show the theorem, we shall see
the explicit crystal structure on 
$\cX:={\mathcal UD}(\chi,T',\theta)$.
Note that ${\mathcal UD}(\chi)=\ZZ^6$ as a set .
Here as  for variables in $\cX$, 
we use the same notations $c,x_0,x_1,\cd,x_5$
as for $\chi$.

For $x=(x_0,x_1,\cd,x_5)\in\cX$, it follows from the 
results in the previous section that the functions
$\wt_i$ and $\vep_i$ ($i=0,1,2$) are given as:
\begin{eqnarray*}
&&\wt_0(x)=2x_0-x_1-x_3-x_5,\,\,
\wt_1(x)=2(x_1+x_3+x_5)-x_0-3x_2-3x_4,\\
&&\wt_2(x)=2(x_2+x_4)-x_1-x_3-x_5.
\end{eqnarray*}
Set
\begin{equation}
\begin{array}{l}
\al\seteq 2x_0+3x_2+x_3,\q
\beta\seteq x_1+3x_2+2x_3+x_5,\q
\gamma\seteq x_0+x_1+3x_3,\,\\
\del\seteq x_0+x_1+x_2+2x_3+x_4,\q
\epsilon\seteq x_0+x_1+2x_2+x_3+2x_4,\q
\\ 
\phi\seteq x_0+3x_2+2x_3,\q
\psi\seteq x_0+x_1+3x_2+3x_4,\q
\xi\seteq x_0+x_1+3x_2+x_3+x_5.
\end{array}
\label{alphabeta}
\end{equation}
Indeed, from the explicit form of $E$ as in the previous
section we have
\[
 {\mathcal UD}(E)=\max(\al,\beta,\gamma,\del,\epsilon,
\phi,\psi,\xi), 
\]
and then
\begin{eqnarray}
&&\vep_0(x)=\max(\al,\beta,\gamma,\del,\epsilon,
\phi,\psi,\xi)-(3x_0+3x_2+x_3),\nn \\
&&\vep_1(x)=\max(x_0-x_1,x_0+3x_2-2x_1-x_3,
x_0+3x_2+3x_4-2x_1-2x_3-x_5),\\
&&\vep_2(x)=\max(x_1-x_2,x_1+x_3-2x_2-x_4).\nn
\end{eqnarray}
Next, we describe the actions of $\eit$ $(i=0,1,2)$.
Set $\Xi_j\seteq{\mathcal UD}(\cC_j)|_{c=1}$ 
($j=1,\cd,5$). Then we have
\begin{eqnarray*}
\Xi_1&=&\max(1+x_0-x_1,x_0+3x_2-2x_1-x_3,
x_0+3x_2+3x_4-2x_1-2x_3-x_5)\\
&&-\max(x_0-x_1,x_0+3x_2-2x_1-x_3,
x_0+3x_2+3x_4-2x_1-2x_3-x_5),\\
\Xi_3&=&\max(1+x_0-x_1,1+x_0+3x_2-2x_1-x_3,
x_0+3x_2+3x_4-2x_1-2x_3-x_5)\\
&&-\max(1+x_0-x_1,x_0+3x_2-2x_1-x_3,
x_0+3x_2+3x_4-2x_1-2x_3-x_5),\\
\Xi_5&=&\max(1+x_0-x_1,1+x_0+3x_2-2x_1-x_3,
1+x_0+3x_2+3x_4-2x_1-2x_3-x_5)\\
&&-\max(1+x_0-x_1,1+x_0+3x_2-2x_1-x_3,
x_0+3x_2+3x_4-2x_1-2x_3-x_5),\\
\Xi_2&=&\max(1+x_1-x_2,x_1+x_3-2x_2-x_4)
-\max(x_1-x_2,x_1+x_3-2x_2-x_4),\\
\Xi_4&=&\max(1+x_1-x_2,1+x_1+x_3-2x_2-x_4)
-\max(1+x_1-x_2,x_1+x_3-2x_2-x_4).
\end{eqnarray*}
Therefore, for $x\in\cX$ we have
\begin{eqnarray*}
&&\til e_1(x)=(x_0,x_1+\Xi_1,x_2,
x_3+\Xi_3,x_4,x_5+\Xi_5),\\
&&\til e_2(x)=(x_0,x_1,x_2+\Xi_2,x_3,x_4+\Xi_4,x_5).
\end{eqnarray*}
We obtain the action $\fit$ ($i=1,2$)
by setting $c=-1$ in ${\mathcal UD}(\cC_i)$.

Finally, we describe the action of $\til e_0$. 
Set
\begin{eqnarray*}
\Psi_0&\seteq&
\max(2+\al,\beta,1+\gamma,1+\del,1+\epsilon,
1+\phi,1+\psi,1+\xi)\\
&&-\max(\al,\beta,\gamma,\del,\epsilon,\phi,\psi,\xi)-1,
\\
\Psi_1&\seteq&
\max(1+\al,\beta,1+\gamma,1+\del,1+\epsilon,
\phi,1+\psi,1+\xi)\\
&&-\max(\al,\beta,\gamma,\del,\epsilon,\phi,\psi,\xi)-1,
\\
\Psi_2&\seteq&
\max(1+\al,\beta,\gamma,1+\del,1+\epsilon,
\phi,1+\psi,1+\xi)\\
&&-\max(\al,\beta,\gamma,\del,\epsilon,\phi,\psi,\xi)-1,
\\
\Psi_3&\seteq&
\max(2+\al,\beta,1+\gamma,1+\del,1+\epsilon,
1+\phi,1+\psi,1+\xi)\\
&&+\max(1+\al,\beta,\gamma,\del,\epsilon,\phi,
\psi,1+\xi)
-\max(1+\al,\beta,\gamma,\del,\epsilon,\phi,
\psi,1+\xi)\\
&&-\max(1+\al,\beta,\gamma,1+\del,1+\epsilon,
\phi,1+\psi,1+\xi)-2,\\
\Psi_4&\seteq&
\max(2+\al,\beta,1+\gamma,1+\del,1+\epsilon,
1+\phi,1+\psi,1+\xi)\\
&&-\max(1+\al,\beta,\gamma,1+\del,1+\epsilon,\phi,
1+\psi,1+\xi)-1,
\\
\Psi_5&\seteq&
\max(2+\al,\beta,1+\gamma,1+\del,1+\epsilon,
1+\phi,1+\psi,1+\xi)\\
&&-\max(1+\al,\beta,\gamma,\del,\epsilon,\phi,
\psi,1+\xi)-1,
\end{eqnarray*}
where $\al,\beta,\cd,\xi$ are as in (\ref{alphabeta}).
Therefore, by the explicit form of $e_0^c$ 
as in the previous section, we have
\begin{equation}
\til e_0(x)=(x_0+\Psi_0,x_1+\Psi_1,x_2+\Psi_2,
x_3+\Psi_3,x_4+\Psi_4,x_5+\Psi_5).
\end{equation}
Now, let us show the theorem.\\
({\sl Proof of Theorem \ref{ultra-d}.})
Define the map 
\[
\begin{array}{cccc}
\Omega\cl&\cX&\longrightarrow& B_\ify,\\
&(x_0,\cd, x_5)&\mapsto& (b_1,b_2,b_3,
\ovl b_3,\ovl b_2,\ovl b_1),
\end{array}
\]
by 
\[
 b_1=x_5,\,\, b_2=x_4-x_5,\,\,
b_3=x_3-2x_4,\,\, \ovl b_3=2x_2-x_3,\,\,
\ovl b_2=x_1-x_2,\,\,
\ovl b_1=x_0-x_1,
\]
and $\Omega^{-1}$ is given by 
\begin{eqnarray*}
&&x_0=b_1+b_2+\frac{b_3+\ovl b_3}{2}+\ovl b_2+\ovl b_1,
\q
x_1=b_1+b_2+\frac{b_3+\ovl b_3}{2}+\ovl b_2,\\
&&
x_2=b_1+b_2+\frac{b_3+\ovl b_3}{2},\,\,
x_3=2b_1+2b_2+b_3,\,\,
x_4=b_1+b_2,\,\,
x_5=b_1,
\end{eqnarray*}
which means that $\Omega$ is bijective.
Here note that 
$\frac{b_3+\ovl b_3}{2}\in \bbZ$ by the 
definition of $B_\ify$.
We shall show that 
$\Omega$ is commutative with actions of $\eit$
and preserves the functions $\wt_i$ and $\vep_i$,
that is, 
\[
 \eit(\Omega(x))=\Omega(\eit x),\q
\wt_i(\Omega(x))=\wt_i(x),\q
\vep_i(\Omega(x))=\vep_i(x)
\q(i=0,1,2).
\]
First, let us check $\wt_i$: Set $b=\Omega(x)$. 
By the explicit forms of $\wt_i$ on 
$\cX$ and $B_\ify$,
we have 
\begin{eqnarray*}
&&\hspace{-10pt}
\wt_0(\Omega(x))=\vp_0(\Omega(x))
-\vep_0(\Omega(x))=2z_1+z_2+z_3+3z_4\q\qq\\
&&=
2(\ovl b_1-b_1)+(\ovl b_2-\ovl b_3)
+(b_3-b_2)+\frac{3}{2}(\ovl b_3-b_3)
=2(\ovl b_1-b_1)+\ovl b_2-b_2
+\frac{\ovl b_3-b_3}{2}\\
&&=2x_0-x_1-x_3-x_5=\wt_0(x),\\
&&\hspace{-10pt}\wt_1(\Omega(x))
=\vp_1(\Omega(x))-\vep_1(\Omega(x))\\
&&=
b_1+(b_3-b_2+({\bar b}_2-{\bar b}_3)_+)_+
-({\bar b}_1+({\bar b}_3
-{\bar b}_2-(b_2-b_3)_+)_+)\\
&&=b_1-\ovl b_1-b_2+\ovl b_2+b_3-\ovl b_3
=2(x_1+x_3+x_5)-x_0-3x_2-3x_4=\wt_1(x),\\
&&\hspace{-10pt}
\wt_2(\Omega(x))
=\vp_2(\Omega(x))-\vep_2(\Omega(x))
=
b_2+\frac{1}{2}({\bar b}_3-b_3)_+-
{\bar b}_2+\frac{1}{2}(b_3-{\bar b}_3)_+\\
&&=b_2-\ovl b_2+\frac{1}{2}(\ovl b_3-b_3)
=2(x_2+x_4)-x_1-x_3-x_5=\wt_2(x).
\end{eqnarray*}
Next, we shall check $\vep_i$:
\begin{eqnarray*}
&&\vep_1(\Omega(x))
={\bar b}_1
+({\bar b}_3-{\bar b}_2+(b_2-b_3)_+)_+\\
&&\,\,=\max({\bar b}_1,{\bar b}_1+
{\bar b}_3-{\bar b}_2,
{\bar b}_1+
{\bar b}_3-{\bar b}_2+b_2-b_3)\\
&&\,\,
=\max(x_0-x_1,x_0+3x_2-2x_1-x_3,
x_0+3x_2+3x_4-2x_1-2x_3-x_5)=\vep_1(x),\\
&&\vep_2(\Omega(x))
={\bar b}_2+\frac{1}{2}(b_3-{\bar b}_3)_+
=\max({\bar b}_2,
{\bar b}_2+\frac{1}{2}(b_3-{\bar b}_3)_+)\\
&&\,\,
\max(x_1-x_2,x_1+x_3-2x_2-x_4)
=\vep_2(x).
\end{eqnarray*}
Before checking $\vep_0(\Omega(x))=\vep_0(x)$,
we see the following formula,
which has been given in \cite[Sect6]{N}.
\begin{lem}
\label{convex}
For $m_1,\cd,m_k\in\RR$ and $t_1,\cd,t_k\in
\RR_{\geq0}$ such that $t_1+\cd t_k=1$,
we have 
\[
 \max\left(m_1,\cd,m_k,\sum_{i=1}^kt_im_i\right)
= \max(m_1,\cd,m_k)
\]
\end{lem}
By the facts
\begin{equation}\label{3tobun}
 \del=\frac{2\gamma+\psi}{3},\q
\epsilon=\frac{\gamma+2\psi}{3},
\end{equation}
and Lemma \ref{convex}, we have
\begin{equation}
\max(\al,\beta,\gamma,\del,
\epsilon,\phi,\psi,\xi)
=\max(\al,\beta,\gamma,\phi,\psi,\xi).
\label{del-ep}
\end{equation}
Here let us see $\vep_0$:
\begin{eqnarray*}
&&\vep_0(\Omega(x))=
-s(b)+\max A-(2z_1+z_2+z_3+3z_4)\\
&&=-x_0+\max(0,z_1,z_1+z_2,z_1+z_2+3z_4,
z_1+z_2+z_3+3z_4,2z_1+z_2+z_3+3z_4)-(\al-\beta)
\\&&
=-x_0+\max
(-2x_0+x_1+x_3+x_5,-x_0+x_3,-x_0+x_1-3x_2+2x_3,
\\&&
\qq\qq \qq\qq-x_0+x_1-x_3+3x_4,-x_0+x_1+x_5,0)\\
&&=
-(3x_0+3x_2+x_3)+\max(x_1+3x_2+2x_3+x_5,
x_0+3x_2+2x_3,x_0+x_1+3x_3,\\
&&\qq x_0+x_1+3x_2+3x_4,
x_0+x_1+3x_2+x_3+x_5,2x_0+3x_2+x_3)\\
&&=-(3x_0+3x_2+x_3)
+\max(\beta,\phi,\gamma,\psi,\xi,\al).
\end{eqnarray*}
On the other hand, we have
\[
\vep_0(x)=-(3x_0+3x_2+x_3)
+\max(\al,\beta,\gamma,\del,\epsilon,
\phi,\psi,\xi).
\]
Then by (\ref{del-ep}), we get 
$\vep_0(\Omega(x))=\vep_0(x)$.

Let us show
 $\eit(\Omega(x))=\Omega(\eit(x))$
($x\in\cX,\,i=0,1,2$).
As for $\til e_1$, set 
\[
 A=x_0-x_1,\,\, B=x_0+3x_2-2x_1-x_3,\,\,
C=x_0+3x_2+3x_4-2x_1-2x_3-x_5.
\]
Then we obtain
$\Xi_1=\max(A+1,B,C)-\max(A,B,C),\,\,
\Xi_3=\max(A+1,B+1,C)-\max(A+1,B,C),\,\,
\Xi_5=\max(A+1,B+1,C+1)-\max(A+1,B+1,C).$
Therefore, we have
\begin{eqnarray*}
&&\Xi_1=1,\,\,\Xi_3=0,\,\,\Xi_5=0,\,\,
\text{if}\,\,A\geq B,C\\
&&\Xi_1=0,\,\,\Xi_3=1,\,\,\Xi_5=0,\,\,
\text{if}\,\,A< B\geq C\\
&&\Xi_1=0,\,\,\Xi_3=0,\,\,\Xi_5=1,\,\,
\text{if}\,\,A, B< C, 
\end{eqnarray*}
which implies
\begin{eqnarray*}
\til e_1(x)=\begin{cases}
(x_0,x_1+1,x_2,\cd,x_5)&\text{if }
A\geq B,C\\
(x_0,\cd,x_3+1,x_4,x_5)&\text{if }
A< B\geq C\\
(x_0,\cd,x_4,x_5+1)&\text{if }
A, B< C
\end{cases}
\end{eqnarray*}
Since $A=\ovl b_1$, 
$B=\ovl b_1+\ovl b_3-\ovl b_2$ and 
$C=\ovl b_1+\ovl b_3-\ovl b_2+b_2-b_3$,
we get ($b=\Omega(x)$)
\begin{eqnarray*}
\Omega({\tilde e}_1(x))=&
\begin{cases}
(\ldots,{\bar b}_2 +1,{\bar b}_1 -1) 
& \text{if ${\bar b}_2 -{\bar b}_3 \geq (b_2 -b_3)_+$}, 
\\  
(\ldots,b_3 +1,{\bar b}_3 -1,\ldots) 
& \text{if ${\bar b}_2 -{\bar b}_3 <0\leq b_3 -b_2$}, 
\\ 
(b_1 +1,b_2 -1,\ldots) 
& \text{if $({\bar b}_2 -{\bar b}_3)_+ <b_2 -b_3$},
\end{cases}
\end{eqnarray*}
which is the same as the action of 
$\til e_1$ on $b=\Omega(x)$ as in Sect.4.
Hence, we have 
$\Omega(\til e_1(x))=\til e_1(\Omega(x))$.

Let us see 
$\Omega(\til e_2(x))=\til e_2(\Omega(x))$.
Set
\[
 L=x_1-x_2, \q M:=x_1+x_3-2x_2-x_4.
\]
Then $\Xi_2=\max(1+L,M)-\max(L,M)$ and 
$\Xi_4=\max(1+L,1+M)-\max(1+L,M)$. Thus,
one has
\begin{eqnarray*}
&&\Xi_2=1,\q\Xi_4=0 \q\text{if }L\geq  M,\\
&&\Xi_2=0,\q\Xi_4=1 \q\text{if }L< M,
\end{eqnarray*}
which means
\[
 \til e_2(x)=\begin{cases}
(x_0,x_1,x_2+1,x_3,x_4,x_5)&\text{if }L\geq M,\\
(x_0,x_1,x_2,x_3,x_4+1,x_5)&\text{if }L< M.
\end{cases}
\]
Since $L-M=x_2-x_3+x_4=\frac{\ovl b_3-b_3}{2}$,
one gets
\[
\Omega( \til e_2(x))=
\begin{cases}
(\ldots,{\bar b}_3 +2,{\bar b}_2 -1,\ldots) 
& \text{if ${\bar b}_3 \geq b_3$}, 
\\
(\ldots,b_2 +1,b_3 -2,\ldots) 
& \text{if ${\bar b}_3 <b_3$},
\end{cases}
\]
where $b=\Omega(x)$. 
This action coincides with the one of
 $\til e_2$ on $b\in B_\ify$ as in Sect4.
Therefore, we get 
$\Omega(\til e_2(x))=\til e_2(\Omega(x))$.

Finally, we shall check 
$\til e_0(\Omega(x))=\Omega(\til e_0(x))$.
For the purpose, we shall estimate the values
$\Psi_0,\cd,\Psi_5$ explicitly.

First, the following cases are investigated:
\begin{eqnarray*}
&({\rm e}1)&\beta>\al, \gamma,\del,\epsilon,
\phi,\psi,\xi,\\
&({\rm e}2)&\beta\leq \phi>\al,\gamma,\del,\epsilon,
\psi,\xi\\
&({\rm e}3)&\beta,\phi\leq \gamma>\al,\del,\epsilon,
\psi,\xi\\
&({\rm e}4)&\beta,\gamma,\del,\epsilon,\phi
\leq\psi>\al,\xi\\
&({\rm e}4')&\beta,\gamma,\epsilon,\phi,\psi
\leq\del>\al,\xi\\
&({\rm e}4'')&\beta,\gamma,\del,\phi,\psi
\leq\epsilon>\al,\xi\\
&({\rm e}5)&\beta,\gamma,\del,\epsilon,
\phi,\psi\leq \xi>\al,\\
&({\rm e}6)&\al\geq \beta,\gamma,\del,\epsilon,
\phi,\psi,\xi.
\end{eqnarray*}
It is  easy to see that each of these
 conditions are equivalent to the conditions 
$(E_1)$-$(E_6)$ in Sect.4, more precisely, we have
$({\rm e}i)\Leftrightarrow\,\, (E_i)$ ($i=1,2,\cd,6$), 
and that (e1)--(e6) cover all cases and they 
have no intersection.
Note that the cases (e4') and (e4'') are
 included in the case (e4) thanks 
to (\ref{3tobun}).

Let us show (e1)$\Leftrightarrow\,(E_1)$:
the condition (e1) means
$\beta-\al=-(2z_1+z_2+z_3+3z_4)>0$, 
$\beta-\gamma=-(z_1+z_2)>0$, 
$\beta-\del=-(z_1+z_2+z_4)>0$, 
$\beta-\epsilon=-(z_1+z_2+2z_4)>0$,
$\beta-\phi=-z_1>0$, $\beta-\psi=
-(z_1+z_2+3z_4)>0$ and 
$\beta-\xi=-(z_1+z_2+z_3+3z_4)>0$, 
which is equivalent to the 
condition $z_1+z_2<0$, $z_1<0$, 
$z_1+z_2+3z_4<0$ and 
$z_1+z_2+z_3+3z_4<0$. This is just 
the condition $(E_1)$.
Other cases are shown similarly.

Under the condition (e1)
($\Leftrightarrow\,(E_1)$),
we have 
\[
 \Psi_0=\Psi_1=\Psi_2=\Psi_4=\Psi_5=-1,\q
\Psi_3=-2,
\]
which means
$\til e_0(x)=(x_0-1,x_1-1,x_2-1,x_3-2,x_4-1,
x_5-1)$. Thus, we have
\[
 \Omega(\til e_0(x))=(b_1-1,b_2,\cd,\ovl b_1),
\]
which coincides with the action of 
$\til e_0$ under $(E_1)$ in Sect.4.
Similarly, we have
\[
\begin{array}{ccc}
{\rm(e2)}&\Rightarrow&
(\Psi_0,\Psi_1,\Psi_2,\Psi_3,\Psi_4,\Psi_5)
=(0,-1,-1,-1,0,0)\\
&\Rightarrow
&\til e_0(x)=(x_0,x_1-1,x_2-1,x_3-1,x_4,x_5),\\
&\Rightarrow&
\Omega(\til e_0(x))
=(b_1,b_2,b_3-1,\ovl b_3-1,\ovl b_2,\ovl b_1+1),
\end{array}
\]
which coincides with the action of $\til e_0$
under $(E_2)$ in Sect.4.
\[
\begin{array}{ccc}
{\rm(e3)}&\Rightarrow&
(\Psi_0,\Psi_1,\Psi_2,\Psi_3,\Psi_4,\Psi_5)
=(0,0,-1,-2,0,0)\\
&\Rightarrow
&\til e_0(x)=(x_0,x_1,x_2-1,x_3-2,x_4,x_5),\\
&\Rightarrow&
\Omega(\til e_0(x))
=(b_1,b_2,b_3-2,\ovl b_3,\ovl b_2+1,\ovl b_1),
\end{array}
\]
which coincides with the action of $\til e_0$
under $(E_3)$ in Sect.4.
\[
\begin{array}{ccc}
\text{(e4)}&\Rightarrow&
(\Psi_0,\Psi_1,\Psi_2,\Psi_3,\Psi_4,\Psi_5)
=(0,0,0,-2,-1,0)\\
&\Rightarrow
&\til e_0(x)=(x_0,x_1,x_2,x_3-2,x_4-1,x_5),\\
&\Rightarrow&
\Omega(\til e_0(x))
=(b_1,b_2-1,b_3,\ovl b_3+2,\ovl b_2,\ovl b_1),
\end{array}
\]
which coincides with the action of $\til e_0$
under $(E_4)$ in Sect.4.
\[
\begin{array}{ccc}
{\rm(e5)}&\Rightarrow&
(\Psi_0,\Psi_1,\Psi_2,\Psi_3,\Psi_4,\Psi_5)
=(0,0,0,-1,-1,-1)\\
&\Rightarrow
&\til e_0(x)=(x_0,x_1,x_2,x_3-1,x_4-1,x_5-1),\\
&\Rightarrow&
\Omega(\til e_0(x))
=(b_1-1,b_2,b_3+1,\ovl b_3+1,\ovl b_2,\ovl b_1),
\end{array}
\]
which coincides with the action of $\til e_0$
under $(E_5)$ in Sect.4.
\[
\begin{array}{ccc}
{\rm(e6)}&\Rightarrow&
(\Psi_0,\Psi_1,\Psi_2,\Psi_3,\Psi_4,\Psi_5)
=(1,0,0,0,0,0)\\
&\Rightarrow
&\til e_0(x)=(x_0+1,x_1,x_2,x_3,x_4,x_5),\\
&\Rightarrow&
\Omega(\til e_0(x))
=(b_1,b_2,b_3,\ovl b_3,\ovl b_2,\ovl b_1+1),
\end{array}
\]
which coincides with the action of $\til e_0$
under $(E_6)$ in Sect.4.
Now, we have 
$\Omega(\til e_0(x))=\til e_0(\Omega(x))$.
Therefore, 
the proof of Theorem \ref{ultra-d} has been
completed.\qed

\bibliographystyle{amsalpha}

\end{document}